\documentclass[12pt]{article}
\usepackage{natbib}
\usepackage{amssymb,amsmath,amsfonts,amssymb}
\usepackage{graphicx}
\usepackage{subcaption}
\usepackage{multicol}
\usepackage{soul}
\input{epsf.sty}
 
\usepackage{comment}

\newcommand{\NI}{\noindent}
\newtheorem{theorem}{\NI{\bf Theorem}}[section]
\newtheorem{lemma}{\NI\bf Lemma}[section]
\newtheorem{prop}{\NI\bf Proposition}[section]
\newtheorem{cor}{\NI\bf Corollary}[section]
\newtheorem{remark}{\NI\bf Remark}[section]
\newtheorem{defn}{\NI\bf Definition}[section]
\newtheorem{example}{\NI\bf Example}[section]
\newcommand{\bt}{\begin{theorem}}
\newcommand{\et}{\end{theorem}}
\newcommand{\bc}{\begin{cor}}
\newcommand{\ec}{\end{cor}}
\newcommand{\bl}{\begin{lemma}}
\newcommand{\el}{\end{lemma}}
\newcommand{\bx}{\begin{example}}
\newcommand{\ex}{\end{example}}
\newcommand{\bea}{\begin{eqnarray}}
\newcommand{\eea}{\end{eqnarray}}
\newcommand{\ben}{\begin{eqnarray*}}
\newcommand{\een}{\end{eqnarray*}}

\newcommand{\be}{\begin{equation}}
\newcommand{\ee}{\end{equation}}

\begin{document}
\begin{center}
{\bf Stochastic optimization  of the Dividend strategy with reinsurance in correlated multiple insurance lines of  business}\\

Khaled Masoumifard$^a$\footnote{Email: k$_{-}$masoumifard@sbu.ac.ir }
and Mohammad Zokaei$^a$\footnote{Corresponding author. Email: zokaei@sbu.ac.ir. Postal Code: 1983969411}
\\
$^a$\,Department of Statistics, Faculty of Mathematical Sciences, Shahid Beheshti University, Tehran, Iran\\
%The authors declare that there is no conflict of interest.
%\today
\end{center}
\begin{abstract}
	The present paper addresses the issue of the stochastic control  of the optimal dynamic reinsurance policy and dynamic dividend strategy, which are state-dependent, for an insurance company that operates under multiple insurance lines of business. The aggregate claims model
	with a thinning-dependence structure is adopted for the risk process.
In the optimization method, the maximum of the cumulative expected discounted dividend payouts with respect to the dividend and reinsurance strategies are considered as value function. 
	This value function  is characterized as the smallest super Viscosity solution  of the associated Hamilton-Jacobi-‌‌‌‌‌‌Bellman (HJB) equation.  
	The finite difference method (FDM) has been utilized for the numerical solution of the value function and the optimal control strategy and the proof for the convergence of this numerical solution to the value function is provided.
	The findings of this paper provide insights for the insurance companies as such that based upon the lines in which they are operating, they can choose a vector of the optimal dynamic reinsurance strategies and consequently transfer  some part of their risks to several reinsurers. The numerical examples in the elicited results show the significance increase in the value function
	 in comparison with the previous findings.
	%in Zhao (2011), Zhao and Balakrishnan (2012),
	%Balakrishnan and Zhao (2013a) and Balakrishnan et al. (2014).
\end{abstract}

%\NI {\bf Keywords:} Cramer-Lundberg process; Optimal reinsurance; Thinning dependence; Hamilton-Jacobi-Bellman equation; Viscosity solution; Dynamic programming principle.

%Generalized Gamma Distribution; Multiple-Outlier
%Scale Model; Hazard Rate Order; Usual Stochastic Order;  Hazard
%Rate Function; Reversed Hazard Rate Function; $P$-Larger Order;
%Parallel System
%%%%%%%%%%%%%%%%%%%%%%%%%%%%%%%%%%%%%%%%%%%%%%%%%%%%%%%%%%%%%%%%%%
%%%%%%%%%%%%%%%%%%%%%%%%%%%%%%%%%%%%%%%%%%%%%%%%%%%%%%%%%%%%%%%%%
%%%%%%%%%%%%%%%%%%%%%%%%%%%%%%%%%%%%%%%%%%%%%%%%%%%%%%%%%%%%%%%%%%
\section{Introduction}
Suppose a insurance company based on a dynamic strategy distributes a ratio of its dividend amongst the shareholders and transfers a part of its risk to a secondary insurance company by a dynamic reinsurance strategy. The dividend and reinsurance strategies are shown as $\{D_t\}_{t\ge 0}$ 
and
$\{R_t\}_{t\ge 0}$, respectively. A paramount issue for an insurance company is the optimization of these strategies. For this reason, first an objective function should be considered and then  $\{D_t\}_{t\ge 0}$ 
and
$\{R_t\}_{t\ge 0}$ strategies should be found as such that the objective function is optimized. A very common function in literature is the cumulative expected discounted dividends which is displayed as $V(.)$.
In the following, we will outline some research on thinning-dependence structure and optimization $V(.)$ with respect to the dividend and reinsurance.\\
		\hl{Optimization of $V(.)$ with respect to the dividend strategy}:
	In 1957, De Finetti (\cite{de1957impostazione}) considered the band strategy for paying dividends to shareholders and addressed the issue of optimizing $V(.)$.
	Gerber	 \cite{gerber1969entscheidungskriterien} demonstrated that if the company's capital is modeled using the Cramer-Landbrug process, then an optimal strategy for paying the dividends is always based on a band strategy, and when the severity of the claim follows an exponential distribution, the band strategy is reduced to the barrier strategy.
Later, the issue of optimizing the dividend distribution strategy has attracted many researchers, for example,
when the company's capital is modeled using Brownian risk,  \cite{gerber1969entscheidungskriterien}, \cite{grandits2007optimal}, and \cite{jeanblanc1995optimization} investigated the problem of optimal dividend distribution, and when the company's capital is modeled using the classical risk process, \cite{zhou2005classical}, \cite{avram2007optimal},  \cite{kyprianou2007distributional} and \cite{loeffen2008optimality}  investigated the barrier model for distributing dividends.
Recently, the issue of stochastic control of dividend strategy by researchers has been investigated in a situation where several insurance companies are co-operating and the company's capital is modeled through multi-dimensional stochastic process( refer to
 \cite{loeffen2008optimality},  \cite{albrecher2017optimal} and \cite{azcue2016optimal}).
	\\
		\hl{Optimization of $ V (.) $ with respect to the dividend strategy and reinsurance strategy}:
	  \cite{azcue2005optimal} 
	 is one of the first articles that optimizes $V(.)$ with respect to the dividend and reinsurance strategies.
	\cite{schmidli2006optimisation}
 has examined two important issues; maximizing $V(.)$ with respect to the dividend strategy and minimizing the ruin probability with respect to the reinsurance strategy.
 \cite{thonhauser2007dividend}  suggested that the optimization of dividend and reinsurance should be based on a value function that depends on the expectation of the dividends and the ruin time. With this approach, when the number of claims follows a Poisson distribution with exponential intensity, an optimal barrier strategy is obtained for both the diffusion model and the Cramer-Lundberg model. Some studies have addressed this issue(e.g.
	  \cite{beveridge2007optimal},   \cite{meng2011optimal} and \cite{zhou2012optimal}).\\
		\hl{ The thinning-dependence structure}:
	An insurance company usually operates in several insurance lines, each of which may behave differently \textit{vis-\`{a}-vis } another. Therefore, it is reasonable for the insurance company to make different decisions about each line.
	 For example,  using proportional reinsurance in one line and excess-of-loss reinsurance in another line or use one type of excess-of-loss in one line and a different type of excess-of-loss insurance in another line. \cite{masoumifard2020optimal}, considering a vector of independent compound Poisson process for the lines, have dealt with optimization of reinsurance strategies corresponding to each line. However, the problem is that usually there is a  correlation between the lines. Therefore a method for modeling this correlation is required. In this regard pertinent models can be found in the related literature.
	 \cite{yuen2002comparing} introduce the idea of claim thinning for investigating this issue and provide a clear explanation of how it works. Some literature on the
	 aggregate claims model
	 with thinning-dependence structure as follows;  \cite{wu2003discrete},
	 \cite{lindskog2003common}, \cite{pfeifer2004modeling}, \cite{bauerle2005multivariate} and \cite{wang2005correlated}. Although research on optimal reinsurance is increasing apace, only a few papers deal with
	 the problem concerning the dependent risks. For two lines of insurance business with common shock dependence, in the dynamic setting,  \cite{bai2013optimal} probed
	 the optimal excess of loss reinsurance to minimize the ruin probability for the diffusion risk model and \cite{liang2016optimal} adopted
	 the variance premium principle to study the optimal proportional reinsurance problem for both
	 the compound Poisson risk model and the diffusion approximation risk model.  For more than two lines
	 of insurance business with the common shock dependence, \cite{yuen2015optimal} considered the objective
	 of maximizing the expected exponential utility and derived the optimal reinsurance
	 strategy not only for the diffusion approximation risk model but also for the compound
	 Poisson risk model. Furthermore, \cite{wei2017optimal}  considered a model where  the claim-number
	 processes amongst the lines of the  insurance business have the thinning-dependence structure. For this risk
	 model, they derived the optimal reinsurance strategies with the objective of maximizing the adjustment coefficient for two commonly-used premium principles.
	 	 
In this paper, the optimization of the control strategy (dividend and reinsurance strategies) is utilized, implementing a dynamic programming approach for an insurance company with several dependent lines. In modeling the risk process, we adopt a thinning dependence structure that covers a large part of the dependency structures  which can be expressed and the control strategy that maximizes the cumulative expected discounted dividends (value function) with the value function have being characterized.

%These papers motivate us to consider the optimal proportional reinsurance with dependent risks
%under the variance premium principle.  By a nonstandard approach, we investigate the conditions
%of existence and uniqueness of the optimal reinsurance strategies. Using the technique of stochastic
%control theory, closed-form expressions for the optimal reinsurance strategy and the value function...
\section{Model and problem formation}
In the classical Cramer-Lundberg process, the reserve $X_t$ of an insurance company can be described by
\begin{eqnarray}
X_t=x+pt-\sum_{i=1}^{N_t}U_i
\end{eqnarray}
where $x\ge 0$ is the deterministic initial capital, $N(t)$ is the number of claims arriving in $[0,\infty]$
 with claim arrival intensity $\beta >0$, and the claim severity $U_i$ are i.i.d random variable with distribution $F$. The premium rate $p$ per unit of time is calculated using the expected value principle with relative safety loading $\eta>0 $; that is, $p=(1+\eta)\beta E(U_i)$. A limitation of the existing in this model is the implicit or explicit assumption that the insurers produce only one type of insurance, even though most insurers produce multiple types 
of coverage (e.g., automobile insurance, general liability insurance, fire insurance, workers' compensation insurance, etc.). 
%%واضح است که نرخ ورود حق بیمه و توزیع مربوط به ادعاها در یک نوع از خط بیمه ممکن است متفاوت از دیگری باشد بعنوان مثال ممکن است در یک نوع از خط بیمه نامه میزان ادعا کم باشد  امافراوانی ادعا زیاد(بیمه اتومبیل و...) و در نوع دیگر میزان ادعا زیاد ولی فراوانی کم باشد(بیمه آتش سوزی و...).  بنابراین  اگر برای هر خط بیمه یک  مدل  کرامر لاندبرگ جداگانه در نظر  بگیریم  و سپس این مدل ها را باهم جمع کنیم   ممکن است فرآیند سرمایه را دقیق تر مدل کنیم. 
The dependency can be introduced between the processes through thinning. Suppose that an insurance company has $n$ $(n \ge 2)$ lines of business and stochastic sources that may cause a claim in at least one of
the $n$ lines are classified into $m$ class. It is assumed that each event in the $k$th
class may cause a claim in the $j$th line with probability $p_{k j}$ for $k = 1, 2, \cdots , m$ and
$j = 1, 2, \cdots , n$.
Consider $N^{(k)}(t)$ the number of events of the $k$th class
occurred up to time $t$ and
$N_j^{(k)} (t)$ the number of claims of the $j$th line up to
time $t$ generated from the events in class $k$.  Then the claim number process of the $j$th line can be written as
$$N_j(t)=\sum_{k=1}^{m}N_j^{(k)} (t)$$
Let us define the process  $X_t=(X_1(t),\cdots,X_n(t))$  such that;
\begin{eqnarray}
X_j(t)=p_jt-\sum_{i=1}^{N_{j}(t)}U_j(i),\qquad \qquad \qquad j=1,\cdots,n
\end{eqnarray}
where the claims severity of the $j$th line $U_j(i)$, $i=1,\cdots,N_j(t)$, are i.i.d random variables with distribution $F_j$.  Let the risk process of the $j$th line of insurance company is modeled by $X_j(t)$.
 The premium rate $p_j$ is calculated using the expected value principle with relative safety loading $\eta>0 $. Given an initial surplus $x$, the surplus $X_t$ of the insurance company at time $t$ can be written as
\begin{eqnarray*}
X(t)=x+\sum_{j=1}^{n}X_j(t)=x+(\sum_{j=1}^{n}p_j)t-\sum_{j=1}^{n}\sum_{i=1}^{N_{j}(t)}U_j(i)
\end{eqnarray*}
To analyze $X(t)$
mathematically tractable,
we assume
\begin{itemize}
	\item[(A1)] 
 the processes $N^{(1)} (t), \cdots , N^{(m)} (t)$ are independent Poisson processes with intensities $\beta_1 , \cdots , \beta_m$, respectively,
\item[(A2)]   $U^{(i)},\,\,i=1,2,\cdots, m,$ are independent random variables where $U^{(i)}$ is the sum of all simultaneous claims that spring from the class $i$ at a moment,
\item[(A3)]  all claims made at different times are independent, and
\item[(A4)]  
for any $z$, $V_z(i)$, $i=1,\cdots,m$, are i.i.d random variable where
$V_z(i)$ is the claim incurred in the line $z$ by class $i$.
\end{itemize}
%For $k\neq k'$, the two vectors of claim-number
%	processes $(N^k (t), N_1^k (t), \cdots , N_n^k (t))$ and $(N^{k'}(t), N_1^{k'}(t), . . . , N_n^{k'}(t))$ are independent.
	%\item[(A2)]  For each $k = 1, \cdots , m,$ $N_1^k (t), \cdots , %N_n^k (t)$ are conditionally independent given
%	$N^k (t)$.
%For more details see  \cite{yuen2002comparing}.
In the following proposition, we show that $X(t)$ is a compound Poisson process.
\begin{prop}
Suppose $S^n=\{1, \cdots,n\}$, and $S_i^n=\{A_{ij}^n; j=1,\cdots,\tiny{\begin{pmatrix}n\\i\end{pmatrix}}  \}$, where $A_{ij}^n$ a subset of a set $S^n$, with exactly $i$ elements. Then, $X(t)$ is still a compound Poisson risk process which can be rewritten as
\begin{eqnarray}\label{classicalmodel}
	X(t)=x+\int_{0}^{t}p_s ds-\sum_{i=1}^{N_t'}Y_i
\end{eqnarray}
where  $N_t'$ is a Poisson process with claim arrival intensity $\beta=\sum_{i=1}^m \beta_i$ and the  $Y_i$ are i.i.d random variable with distribution 
$$G(\alpha)=\sum_{j=1}^{n}\sum_{k=1}^{\tiny{\begin{pmatrix}n\\j\end{pmatrix}}}\sum_{i=1}^{m}\frac{\beta_i}{\sum_{i=1}^{m}\beta_i}\prod_{z\in A_{jk}^n}p_{iz}\prod_{z\in S^n-A_{jk}^n}(1-p_{iz}) P(\sum_{z \in A_{jk}^n}V_z\le \alpha) $$
where
$V_z$ is the claim incurred in the line $z$. 
\end{prop}
{\bf  Proof }:\,\, The sum of claims arrived to the insurance company that spring from the class  $i$ up to
time $t$ will be displayed as $$Z^{(i)}(t)=\sum_{k=1}^{N^{(i)}(t)} Z_k^{(i)},$$
where $Z_k^{(i)}$ is the sum of all simultaneous claims that spring from the class $i$ at a moment. 
By using assumptions (A1) and (A3), it is quit axiomatic $Z^{(i)}(t)$ is a compound Poisson process where $N^{(i)}$ is a Poisson distribution with intensity parameter  $\beta_i$ and has the following distribution:
 $$F_i(\alpha)=\sum_{j=1}^{n}\sum_{k=1}^{\tiny{\begin{pmatrix}n\\j\end{pmatrix}}}\prod_{z\in A_{jk}^n}p_{iz}\prod_{z\in S^n-A_{jk}^n}(1-p_{iz})P(\sum_{z \in A_{jk}^n}V_z(i)\le \alpha),$$
 wherein $P(\sum_{z \in A_{jk}^n}V_z(i)\le \alpha)$ indicates the probability of occurrences which arrive from class $i$, having resulted in the claims at the set $A_{jk}^n $ of the lines. From assumption (A4), $P(\sum_{z \in A_{jk}^n}V_z(i)\le \alpha)=P(\sum_{z \in A_{jk}^n}V_z\le \alpha)$. Given the fact that $Z^{(1)}(t),\cdots, Z^{(m)}(t)$ are independent (assumptions (A1) and (A2)), $\sum_{i=1}^mZ^{(i)}(t)$ is also a compound Poisson process like $\sum_{i=1}^{N(t)} Y_i$, where $N(t)=\sum_{i=1}^{m}N_i(t)$ is a Poisson distribution with intensity parameter $\sum_{i=1}^m \beta_i$ and $Y$ has the distribution function $\sum_{i=1}^{m}\frac{\beta_i}{\sum_{i=1}^m \beta_i} F_i$ (refer to  \cite{klugman2012loss}, Theorem 7.5).  $\square$\\
 
Note that the risk process $X(t)$, defined in (\ref{classicalmodel}), is an  càdlàg stochastic
process and satisfies the Markov property. We can describe this model  by defining  $(\Omega, \Sigma, ({\cal{F}}_t)_{t\ge 0}, P)$ as the smallest filtered probability space produced by $X(t)$ (refer to section 1.1 of \cite{azcue2014stochastic}).
\begin{remark}\label{remaek_example}
	If $m = n$ and $p_{kk} = 1$ for $k = 1, \cdots , n$, then $X(t)$ of (\ref{classicalmodel}) is the
	risk model considered by \cite{yuen2002comparing}. For example, let $k=3$:
	\begin{align*}
	G(\alpha)=&
	 \frac{1}{\sum_{i=1}^{3} \beta_i}\bigg(\beta_1(1-p_{12})(1-p_{13})P(U_1\le \alpha)+\beta_2(1-p_{21})(1-p_{23})P(U_2\le \alpha)\\
	&+{\beta_3} (1-p_{31})(1-p_{32})P(U_3\le \alpha)+
	\big(\beta_1p_{12}(1-p_{13})+\beta_2p_{21}(1-p_{23})\big)P(U_1+U_2\le \alpha)\\
	&+
	\big(\beta_1p_{13}(1-p_{12})+\beta_2p_{31}(1-p_{32})\big)P(U_1+U_3\le \alpha)\\
	&+
	\big(\beta_2p_{23}(1-p_{21})+\beta_3p_{32}(1-p_{31})\big)P(U_2+U_3\le \alpha)\\
	&+\big(\beta_1p_{11}p_{12}p_{13}+\beta_2p_{21}p_{22}p_{23}+\beta_3 p_{31}p_{32}p_{33}\big)P(U_1+U_2+U_3\le \alpha)\bigg).
	\end{align*}
	If $n = 3, m =7$, $p_{12} = p_{21} = p_{13}=p_{31}=p_{23}=p_{32}=  0$, $p_{41} = p_{42} = 1$, $p_{52} = p_{53} = 1$, $p_{61} = p_{63} = 1$, $p_{71} = p_{72}  =p_{73}=1$and $p_{11} = p_{22} = p_{33}=1$, then $X (t)$ of (\ref{classicalmodel}) is
	the risk model with common shock for three dependent lines of business;
	\begin{align*}
		G(\alpha)=\frac{1}{\sum_{i=1}^{7}\beta_i}\bigg(&\beta_1P(U_1\le \alpha)+\beta_2 P(U_2\le \alpha)+\beta_3P(U_3\le \alpha)+\beta_4P(U_1+U_2\le \alpha)\\
		&+\beta_5P(U_2+U_3\le \alpha)+\beta_6P(U_1+U_3\le \alpha)+\beta_7P(U_1+U_2+U_3\le \alpha)\bigg)
		\end{align*}
	For $n > 3$, more general risk models with common shock can also
	be constructed from (\ref{classicalmodel}) by choosing the values of $m$ and $p_{k j}$ appropriately.
	
\end{remark}

\subsection{Control strategy}
A control strategy is a process $\pi=(\boldsymbol{R}, D)$ where $\boldsymbol{R}$ is a vector of  reinsurance strategies and  $D_t$ is a dividend strategy.
Reinsurance can be an effective way to manage risk by transferring risk from an
insurer to a second insurer (referred to as
the reinsurer).
A reinsurance contract is an agreement between an insurer and a reinsurer under which, claims that arise are shared between the insurer and reinsurer.
Let a Borel measurable function $R:[0,\infty)\longrightarrow[0,\infty)$,  be called retained loss function, describing the part of the claim that the company pays and satisfies $0\le R(\alpha)\le \alpha$. The reinsurance company covers $\alpha - R(\alpha)$, where the severity of the claim is $\alpha$.  Now,  to reduce the risk exposure of the portfolio, assume that the insurer can take reinsurances in a dynamic way for some insurance lines, each of these reinsurances is indexed by $\{1,\cdots,n\}$.  We denote by $\boldsymbol{\mathcal{F}}$ the vector $(\mathcal{F}_1,\cdots,\mathcal{F}_n)$, in which $\mathcal{F}_i$ is the family of retained loss functions associated to the reinsurance policy in $i$'th line. 
% همچنین با توجه به متفاوت بودن فرآیندمربوط به خط بیمه ها ممکن است قراددادهای متفاوتی برای هر خط بیمه در نظر گرفته شودو یا برای تعدادی از خط بیمه ها از بیمه اتکایی استفاده نکنیم و تمام تعهدات را شکت بیمه به تنهایی ب\زیرد. 
Thus, the reinsurances control strategy is a collection 
$\boldsymbol{R}=(\boldsymbol{R}_t)_{t \ge 0}=({R_1}_t,\cdots, {R_n}_t)_{t \ge 0}$
of the {\Large }vector functions
$\boldsymbol{R}_t:\boldsymbol{\Omega}\rightarrow \boldsymbol{\mathcal{F}}$  for any $t\ge 0$. 

 Well-known reinsurance types are:
\begin{itemize}
	\item[(1)] Proportional reinsurance with $R_P(\alpha)=b\alpha$.
	\item[(2)] Excess of loss reinsurance (XL) with $R_{XL}(\alpha,M) = \min \{\alpha, M \}$, $0\le M \le \infty$.
	\item[(3)]Limited XL reinsurance (LXL) with $R_{LXL}(\alpha,M)=\min\{\alpha,M\}+(\alpha-M-L)^+$, $0\le M,L \le \infty$.
\end{itemize}
The numbers $M$ and $L$ are named priority and limit, respectively. 

A dividend strategy is a process $D=(D_t)_{t\ge 0}$ where $D_t$ is the cumulative amount of dividends paid out by the reinsurance.
Denote by 
$\Pi_{x}$ the set of all control strategies with initial surplus $x \ge 0$.
Now, for any $\pi \in \Pi_{x}$, the surplus process can be written as
\begin{eqnarray}\label{2.3f}
X^{\pi}(t)&=x+\sum_{j=1}^{n}\left(\int_{0}^{t}p_{{R_j}_s}ds-\sum_{i=1}^{N_j(t)}{R_j}_{\tau_i}(U_j(t))\right)-D_t\nonumber\\
&=x+\int_{0}^{t}\sum_{j=1}^{n} p_{{R_j}_s}ds-\sum_{j=1}^{m}\sum_{i=1}^{N_j(t)}{R_j}_{\tau_i}(U_j(i))-D_t.
\end{eqnarray}
It is easy to see that $X^{\pi}(t)$ is equivalent to the following process
\begin{eqnarray}\label{2.3'}
x+\int_{0}^{t}p_{{\boldsymbol{R}}_s}ds-\sum_{i=1}^{N_t'}Z_i-D_t,
\end{eqnarray}
%\begin{eqnarray}
%Y_t\stackrel{st}{=}x+\sum_{k=1}^{n}\big(\int_{0}^{t}[p_k-\lambda_k\pi_k(U^{(k)}-\mathcal{F}_k(U^{(k)},u_s^k))]ds-\sum_{i=1}^{N_t^%{(k)}}Z_i
%\end{eqnarray}
where  $N_t'$ is a Poisson process with claim arrival intensity $\beta=\sum_{i=1}^m \beta_i$ and the  $Z_i$ are i.i.d random variable with distribution 
\begin{align}
\label{Gdistribution}
G^{\boldsymbol{R}}(\alpha)=\sum_{j=1}^{n}\sum_{k=1}^{\tiny{\begin{pmatrix}n\\j\end{pmatrix}}}\bigg[\sum_{i=1}^{m}\frac{\beta_i}{\sum_{i=1}^{m}\beta_i}\prod_{z\in A_{jk}^n}p_{iz}\prod_{z\in S^n-A_{jk}^n}(1-p_{iz})\bigg]F_{{\boldsymbol{R}}_{A_{jk}^n}}(\alpha),\end{align}
where $F_{{\boldsymbol{R}}_{A_{jk}^n}}(\alpha)=P(\sum_{z \in A_{jk}^n}R_z(U_z)\le \alpha) $.
%\begin{eqnarray}
%Y_t=x+\sum_{k=1}^{n}\big(\int_{0}^{t}[p_k-\lambda_k\pi_k(U^{(k)}-\mathcal{F}_k(U^{(k)},u_s^k))]ds-\sum_{i=1}^{N_t^{(k)}}\mathcal{%R}_k(U_i^{(k)},u_{\tau_i}^k)\big)
%\end{eqnarray}
%\begin{eqnarray}
%Y_t=x+\sum_{k=1}^{n_1}(p_kt-\sum_{i=1}^{N_t^{(k)}}U_i^{(k)})+\sum_{k=n_1+1}^{n}\big(\int_{0}^{t}[p_k-\lambda_k\pi_k(U^{(k)}-\mathc%al{R}_k(U^{(k)},u_s^k))]ds-\sum_{i=1}^{N_t^{(k)}}\mathcal{F}(U_i^{(k)},u_{\tau_i}^k)\big)
%\end{eqnarray}
The time of ruin for this process is defined by 
\begin{align}
\label{2.4}
\tau^{\pi}=\inf\left\{t \ge 0: X^{\pi}(t)<0 \right\}.
\end{align}
In this paper, we assume that the reinsurance calculates its premium using the expected value principle with reinsurance safety loading factor $\eta_1 \ge \eta >0$:
$$q_R=(1+\eta_1)(\sum_{i=1}^{m}\beta_i)E(Y-Z)=(1+\eta_1)(\sum_{i=1}^{m}\beta_i)\bigg(\int_{0}^{\infty}\alpha d G(\alpha)-\int_{0}^{\infty}\alpha d G^{\boldsymbol{R}}(\alpha)\bigg),$$
and so $p_R=p-q_R$, where   $$p=(1+\eta)(\sum_{i=1}^{m}\beta_i)E(Y)=(1+\eta)(\sum_{i=1}^{m}\beta_i)\int_{0}^{\infty}\alpha d G(\alpha).$$
Let  $(\Omega, \Sigma, ({\cal{F}}_t)_{t\ge 0}, P)$ is the smallest filtered probability space produced by $X(t)$. The control strategy $\pi=(\boldsymbol{R}, D)$ is admissible if satisfies 
\begin{itemize}
	\item 
the process $\boldsymbol{R}=(\boldsymbol{R}_t)_{t\ge 0}$ is predictable, that is, the function $(\omega, \alpha, t) \rightarrow \boldsymbol{R}_t (\omega, \alpha)$ is $\Sigma \times Borel \times Borel$    		measurable and the function
$\omega \rightarrow \boldsymbol{R}_t (\omega, \alpha)$  is ${\cal{F}}_{t^{-}}$ measurable for every $t\ge 0$ and  $\alpha \ge 0$, and
\item  the process $D=(D_t)_{t \ge 0}$ is predictable, nondecreasing and càglàd (left continuous with right limits).

\end{itemize}
\subsection{Problem formation}
Regarding an admissible control strategy  $(\pi_t)_{t\ge 0}$ and an initial reserve $x\ge 0$, we define the following  value function:
	\begin{align*}
	V^{\pi}(x)=E_x\bigg[\int_{0}^{\tau^{\pi}}e^{-\delta s}dD_s \bigg].
	\end{align*}

In the case that insurance company considers only one type of reinsurance contract for all the risks which it encounters and insurance portfolio was modeled by Cramer-Lundberg model; \cite{azcue2005optimal} was considering  $V^{\pi}(x)$, and found a general dynamic control strategy that maximizes this value function. Our aim in this paper is to extend this result for the model described earlier, in other words, we are looking for 
\begin{eqnarray}\label{valuefunction}
V(x)=\underset{\boldsymbol{\pi}\in \Pi_x}{\sup}\,\,V^{\pi}(x).
\end{eqnarray}
The method of  dynamic programming is used to characterize the optimal value function (\ref{valuefunction})
 and the  corresponding  optimal reinsurance strategies via the Hamilton-Jacobi-Bellman (HJB) equation. To obtain the Hamilton-Jacobi-Bellman (HJB) equation associated with the value function (\ref{valuefunction}), we need to state the so-called Dynamic Programming Principle (DPP). Similar to the Proposition 3.1 of \cite{azcue2005optimal}, we can obtain the following result.
\begin{lemma}\label{DPPVL}Given any initial $x \ge 0$, we have 
	\begin{align}\label{DPPV}
	V(x)=\sup_{\pi\in \Pi_x}E_x\bigg[\int_{0}^{\tau^{\pi}\wedge h}e^{-\delta t}dD_t+e^{-\delta(\tau^{\pi}\wedge h)}V\big(X^{\pi}(\tau^{\pi}\wedge h)\big)\bigg].
	\end{align}
\end{lemma}
We  now  deduce  the   HJB equation  assuming  some  regularity  on  $V$. For any continuously differentiable function $u$ defined in $R_+$, we define the discounted infinitesimal generator $\mathsf{g}$  of  the controlled process 
$X^{\pi}(t)$ by
\begin{align*}
\mathsf{g}( X^{\pi}(\tau^{\pi}\wedge t),u)(x)=lim_{t\longrightarrow 0}\frac{E_x(e^{-\delta t}u(X^{\pi}(\tau^{\pi}\wedge t)))-u(x)}{t}.
\end{align*}
Assume that $V$ is continuously differentiable at x. Given any $d\ge  0$ and any $\boldsymbol{R} \in \boldsymbol{\mathcal{F}}$, let us consider the admissible control strategy $((D_t)_{t\ge 0}, (\boldsymbol{R})_{t \ge 0})$.  Now, similarly as  in section 1.4 of \cite{azcue2014stochastic}
\begin{align*}
\mathsf{g}( X^{\pi}(\tau^{\pi}\wedge t),u)(x)=&\mathsf{g}( X^{\pi}(\tau_1^{\pi}\wedge t),u)(x)\\
=&(p_{\boldsymbol{R}}-d)V'(x)-(c+\sum_{i=1}^m \beta_i)V(x)+(\sum_{i=1}^{m}\beta_i)\int_{0}^{x}V(x-\alpha)dG^{\boldsymbol{R}}(\alpha)
\end{align*}
where $\tau_1^{\pi}$ denotes the first claim arrival.
 Note that from Lemma \ref{DPPVL} we have:
\begin{align*}
0\ge E_x\bigg[\int_{0}^{\tau^{\pi}\wedge t}e^{-\delta s}dD_s+e^{-\delta(\tau^{\pi}\wedge t)}V\big(X^{\pi}(\tau^{\pi}\wedge t)\big)\bigg]-V(x).
\end{align*}
Hence, dividing the above inequality by $t$ and taking $t\searrow 0$ gives
\begin{eqnarray*}
	\underset{d\in\mathbb{R}_+,\boldsymbol{R}\in  \boldsymbol{\mathcal{F} }}{sup}\,\,\left(d(1-V'(x))+{{\cal{L}}}_{\boldsymbol{\mathcal{F} }}(V)(x)\right)\le 0,
\end{eqnarray*}
where 
\begin{align}\label{LV} {{\cal{L}}}_{\boldsymbol{R}}(V)(x)=p_{\boldsymbol{R}}V'(x)-(\delta+\sum_{i=1}^m \beta_i)V(x)+(\sum_{i=1}^m \beta_i)\int_{0}^{x}V(x-\alpha)dG^{\boldsymbol{R}}(\alpha).\end{align}
So, the HJB equation can be written as
\begin{eqnarray}\label{hjbeqV}
\max\{1-V'(x),\underset{\boldsymbol{\mathcal{F} }}{sup}\,\,{{\cal{L}}}_{\boldsymbol{R}}(V)(x)\}=0.
\end{eqnarray}
\subsection{Dividend band strategy with reinsurance}
Let ${\cal{A}}$, ${\cal{B}}$, and ${\cal{C}}$ are disjoint sets with ${\cal{A}}\cup {\cal{B}}\cup {\cal{C}}=\mathbb{R}_+$, we say ${\cal{P}}=({\cal{A}},{\cal{B}},{\cal{C}})$ is a band partition
 if ${\cal{A}}$ is closed, bounded, and nonempty;  ${\cal{C}}$ is open from the right; ${\cal{B}}$ is open from the left, the lower limit of any connected component of ${\cal{B}}$ belongs to ${\cal{A}}$, and there exists $b\ge 0$ such that $(b,\infty)\in {\cal{B}}$.
\begin{defn}\label{defi5.3}
	Consider an initial surplus $x\ge 0$, a stationary reinsurance control $\boldsymbol{r}^x=(r_1^x,\cdots, r_n^x)$, and a band partition ${\cal{P}}=({\cal{A}},{\cal{B}},{\cal{C}})$. An admissible control strategy $\pi^x=(\boldsymbol{R}^x, D^x)=(\boldsymbol{R}_t^x, D_t^x)_{t \ge 0}\in \Pi_{x}^{\pi}$ is  defined as follows,
	\begin{itemize}
		\item if $x \in {\cal{A}}$,  we set $D_t^x=p_{\boldsymbol{r}}t= \sum_{i=1}^{n}p_{r_i^x}t$ and $\boldsymbol{R}_t^x=\boldsymbol{r}^x$. Afterward, follow the strategy corresponding to initial surplus $x- \boldsymbol{r}^x(U_1)$ where $U_1$ is the severity of first claim and $\boldsymbol{r}^x(U_1)=\sum_{i=1}^{n}r_i^x(U_1)I_{U_1\in L_i}$,
		where ${U_1\in L_i}$ indicates that $U_1$ is a claim from the line $i$,
		\item if $x \in {\cal{B}}$,  there exists $x_0 \in {\cal{A}}$ such that  $(x_0,x) \subset {\cal{B}}$, then we set $L_0^x=x-x_0$ and $\boldsymbol{R}_0^x=\boldsymbol{r}^x$. Afterward, follow the strategy corresponding initial surplus $x_0$ and
		\item if $x \in {\cal{C}}$,  there exists $x_1 \in {\cal{A}}$ such that  $(x_1,x) \subset {\cal{B}}$. Then $D_t^x=0$ and $\boldsymbol{R}_t^x=\boldsymbol{r}^{X_{t^-}}$ up to $\tau'=\inf\{t: X_t\notin {\cal{C}}\}$. Afterward, follow the strategy corresponding initial surplus $X_{\tau'}$.
	\end{itemize}
The family $\pi({\cal{P}},\boldsymbol{r})=\{(\boldsymbol{R}^x, D^x)\in \Pi_{x}^{\pi}, x\ge 0\}$ is called the reinsurance band strategy associated with ${\cal{P}}$ and $\boldsymbol{r}$.
\end{defn}
We need to define the following operator, 
\begin{align}\label{Lambda} \Lambda_{\boldsymbol{R}}(W)(x)=p_{\boldsymbol{R}}-(c+\sum_{i=1}^m \beta_i)W(x)+(\sum_{i=1}^m \beta_i)\int_{0}^{x}W(x-\alpha)dG^{\boldsymbol{R}}(\alpha).\end{align}
It is easy to see that, if $V'(x)=1$  then 
${{\cal{L}}}_{\boldsymbol{R}}(V)(x)=\Lambda_{\boldsymbol{R}}(V)(x)$. Consider
a stationary reinsurance control $\boldsymbol{r}^x=(r_1^x,\cdots, r_n^x)$  and a band partition ${\cal{P}}=({\cal{A}},{\cal{B}},{\cal{C}})$. In following proposition, using two operators ${{\cal{L}}}_{\boldsymbol{R}}$  and $\Lambda_{\boldsymbol{R}}$,  a verification result is given for the value function $V_{{\cal{P}},\boldsymbol{r}}$ of the reinsurance band strategy $\pi({\cal{P}},\boldsymbol{r})$.
\begin{prop}\label{prop5.5}
Let $W$ is left continuous at the upper limit of the connected component of ${\cal{C}}$, right continuous at the lower limits of the connected components of ${\cal{B}}$, has derivative equal to 1 on ${\cal{B}}$, and an almost-everywhere solution ${{\cal{L}}}_{\boldsymbol{r}^x}(W)=0$  in the connected components of ${\cal{C}}$, and a solution of $\Lambda_{\boldsymbol{r}^x}(W)=0$ in ${\cal{A}}$. Then $W$ is equivalent to $V_{{\cal{P}},\boldsymbol{r}}$.
\end{prop}
\subsection{Viscosity solutions}
The dynamic programming method is a cogent means
to scrutinize the stochastic control problems through  the HJB  equation.
In (\ref{hjbeqV}),  we have obtained the associated equation to the
value function (\ref{valuefunction}). Nonetheless, in the classical approach, this method is adopted only when it is assumed \textit{a priori} 
that optimal value functions are smooth enough. In general,  the optimal value function is not expected  to be smooth enough to satisfy these equations in the classical sense. 
These call for the desideratum of week notation of solution of the HJB equation: the theory of viscosity solutions. Let us define this notion(see \cite{azcue2014stochastic}).
 \begin{defn}\label{Viscositydefi}
 	Let $Z$ be the set of locally Lipschitz functions in $\mathbb{R}_+$. Given a function  $L(x_1,x_2,x_3,g):\mathbb{R}^3\times Z\rightarrow \mathbb{R}$ and a domain $J\subset \mathbb{R}_{+}$, consider the first-order differential equations of the form 
 	\begin{align}
 	\label{3.6}
 L(x,u(x), u'(x),u)= 0,\,\, with\,\, x \in \mathbb{R}_+.
 	\end{align}
 	A function $\underline{u}:J\rightarrow \boldsymbol{R}$ is a viscosity subsolution of the differential equation (\ref{3.6}) at $x\in J$ if $\underline{u}$ is locally Lipschitz and $L(x,\underline{u}(x), \bar{d},\underline{u})\ge 0$ for all $\bar{d}\in D^+(\underline{u})(x)$, where $ D^+(\underline{u})(x)$ is the set of all the super-differentials, that is,
 	$$\liminf_{h\rightarrow 0^{-}}\frac{u(x+h)-u(x)}{h}\ge \bar{d}\ge \limsup_{h\rightarrow 0^{+}}\frac{u(x+h)-u(x)}{h}.$$
 		A function $\underline{u}:J\rightarrow \boldsymbol{R}$ is a viscosity supersolution of the differential equation (\ref{3.6}) at $x\in J$ if $\bar{u}$ is locally Lipschitz and $L(x,\bar{u}(x), \underline{d},\bar{u})\ge 0$ for all $\underline{d}\in D^-(\underline{u})(x)$, where $ D^-(\underline{u})(x)$ is the set of all the sup-differentials,that is,
 	$$\liminf_{h\rightarrow 0^{+}}\frac{u(x+h)-u(x)}{h}\ge \underline{d}\ge \limsup_{h\rightarrow 0^{-}}\frac{u(x+h)-u(x)}{h}.$$ Finally, a function $u:J\rightarrow \boldsymbol{R}$ is a viscosity solution of  (\ref{3.6}) at $x \in J$, if it is both viscosity subsolution and supersolution.
 \end{defn}
There is an equivalent formulation for viscosity subsolution and supersolution.
\begin{defn}
	We say that a  function $\underline{u}:[0,\infty)\rightarrow \mathbb{R} $ is a viscosity subsolution of (\ref{3.6}) at $x \in (0,\infty)$ if it is locally Lipschitz and any continuously differentiable function $\varphi:(0,\infty)\rightarrow \mathbb{R} $ (called test function), with $\varphi(x)=\underline{u}(x)$ such that $\underline{u}-\varphi$ reaches the maximum at $x$, satisfies 
$L(x,\varphi(x), \varphi'(x),\varphi)\ge 0$. 	We say that a continuous function $\bar{u}:[0,\infty)\rightarrow \mathbb{R}$ is a viscosity supersolution of (\ref{3.6}) at $x \in (0,\infty)$ if it is locally Lipschitz  any continuously differentiable function $\phi:(0,\infty)\rightarrow \mathbb{R}$ (called test function),  with $\phi(x)=\bar{u}(x)$ such that $\bar{u}-\phi$ reaches the maximum at $x$ satisfies 
	$L(x,\phi(x), \phi'(x),\phi)\le 0$.
\end{defn}

\section{Main results}
In this section, we state a comparison result between viscosity subsolutions
and supersolutions of (\ref{hjbeqV})  with a suitable boundary condition that gives us the uniqueness of viscosity solution. Also, we characterize the optimal value function as the smallest supersolution of the HJB equation.
Before stating the main results, we need the following lemmas. The proofs of these lemmas are very similar to the same results in section 2.1.2 of  \cite{azcue2005optimal},  we thus omit them.
 \begin{lemma}\label{prop1ofcani}
	For $x \ge 0$, the optimal value function $V(x)$ is well defined and admits the following bound:
	$$x+ \frac{(1+\eta)\mu}{\sum_{i=1}^m\beta_i+\delta}\sum_{i=1}^n \beta_i \le V(x)\le x+ \frac{(1+\eta)\mu}{\delta}\sum_{i=1}^n \beta_i.
	$$
\end{lemma}
\begin{lemma}\label{lemma1ofcani}
	The optimal value function $V$ is increasing and locally Lipschitz in $[0,\infty)$ and for $y>x \ge 0$ satisfies
	$$y-x \le V(y)-V(x)\le (\sum_{i=1}^{n}\beta_i)\frac{V(x)}{\sum_{i=1}^{n}p_i}(y-x).$$
\end{lemma}
Since the optimal value function $V$ is locally Lipschitz but possibly not differentiable at some points, we cannot say that $V$ is a solution of the HJB equation, we prove instead that
$V$ is a viscosity solution of the corresponding HJB equation.
\begin{theorem}\label{V-viscosity}
	$V$ is a viscosity solution of the HJB equation (\ref{hjbeqV}) at any $x \in \mathbb{R}_+$.
\end{theorem}
Now, we first provide  the  comparison principle for viscosity solution of the HJB equation (\ref{hjbeqV}). 
This result implies the uniqueness among a certain class $(L)$ of the viscosity solution of (\ref{hjbeqV}) to which the optimal value function belongs. To be more precise, we introduce the following definition.
\begin{defn}\label{growthA1}
	We say that a function $u:[0,\infty)\rightarrow \mathbb{R}$  belongs class $L$ if satisfies
	\begin{itemize}
		\item[(i)] $u$ is locally Lipschitz,
		\item[(ii)] if $0\le x<y$, then $u(y)-u(x)\ge y-x$, and
		\item[(iii)]   there exists a constant $k > 0$ such that $u(x) \le x+ k$ for all  $x \in [0,\infty)$.
	\end{itemize}
	We also define;
	$L^*=$\{ $u:\,\, u$ is viscosity solution of (\ref{hjbeqV}) and belongs to $L$\}.
\end{defn}
It is interesting to note that if $u$ is of class $L$, then $u$ is strictly positive, linearly bounded, nondecreasing and absolutely continuous.  Absolute continuity follows from the local Lipschitz continuity on a compact set.
Clearly, by Lemma \ref{prop1ofcani}   and \ref{lemma1ofcani} the optimal value function $V$  belongs to $L$.
\begin{prop}\label{prop4.3.}
	(Comparison principle) Let $\underline{u}$ and $\bar{u}$ are the sup and super-viscosity solution of (\ref{hjbeqV}) respectively. If both $\underline{u}$ and $\bar{u}$ are of class $L$, then  $\underline{u}\le \bar{u}$ for $x=0$ implies that $\underline{u}\le \bar{u}$ in $\mathbb{R}_+$.
\end{prop}
The comparison principle states that there is at most one viscosity solution of (\ref{hjbeqV}) with boundary condition at zero among all the functions in $L$ with the same boundary condition.
Since if $K(x)$ and $V(x)$ are two viscosity solutions of (\ref{hjbeqV}) with $V(0)=K(0)$, then $V$ is a viscosity sub-solution and $K$ is a viscosity super-solution, then according to the above proposition $V \le K$. Also, with a similar argument, $V \ge K$ and as a result  $V(x)=K(x)$.
Therefore, by having this result, if we know $V(0)$ then $V(x)$ would be characterized. 
However,
 the problem here is that $V(0)$ is not known a priori; therefore, this result is not enough to characterize $V(x)$. 
By the following proposition, we can finalize the characterization $V$. 
\begin{prop}\label{prop4.4}
	The optimal value function $V(x)=\underset{\boldsymbol{\pi}\in \Pi_x}{\sup}\,\,V^{\pi}(x)$ is the smallest viscosity supersolution of (\ref{hjbeqV}) that belongs to $L$.
\end{prop}
These results allows us to characterize V as the unique viscosity solution of (\ref{hjbeqV}) with boundary condition
$V(0)=\inf_{u\in L^*}u(0)$.
From the previous proposition we can deduce the usual viscosity verification
result:
If we can find a stationary reinsurance  strategy $\pi=(\boldsymbol{R}^x, D^x)\in \Pi_{x}$
such that $V^{\pi}$ is a viscosity supersolution of (\ref{hjbeqV}), then $V(x)=V^{\pi}(x)$;
because $V(x)\ge V^{\pi}(x)$ and by above proposition  $V(x)$ is the smallest viscosity supersolution of (\ref{hjbeqV}). Now we can show that the optimal control strategy is a reinsurance band strategy.

 \begin{prop}\label{prop5.7}	Let the vector $\boldsymbol{\mathcal{F}}=(\mathcal{F}_1,\cdots,\mathcal{F}_n)$, where $\mathcal{F}_i$ is one of the reinsurance families $\mathcal{F}_P$, $\mathcal{F}_{XL}$ and $\mathcal{F}_{LXL}$. Then,
 	${\cal{P}}^*=({\cal{A}}^*,{\cal{B}}^*,{\cal{C}}^*)$ is a band partition, where 
 	\begin{align*}
 	{\cal{A}}^*&=\{x\in \mathbb{R}_+ \,\,\text{suct that}\,\,\sup_{\boldsymbol{{R}} \in \boldsymbol{\mathcal{F}}}\,\,\Lambda_{\boldsymbol{R}}({V})(x)=0  \},
 	\text{ and}\\
 	{\cal{B}}^*&=\{x \in (0,\infty) \,\,\text{suct that}\,\, V'(x)=1\,\,\text{and}\,\,sup_{\boldsymbol{{R}} \in \boldsymbol{\mathcal{F}}}\,\,\Lambda_{\boldsymbol{R}}({V})(x)<0\}\,\,\text{and}\,\,{\cal{C}}^*=({\cal{A}}^*\cup{\cal{B}}^*)^c .\end{align*}
 \end{prop} 
In the rest of this paper, whenever ${\cal{P}}^*=({\cal{A}}^*,{\cal{B}}^*,{\cal{C}}^*)$  is used, it refers to the ${\cal{P}}^*=({\cal{A}}^*,{\cal{B}}^*,{\cal{C}}^*)$ introduced in the previous proposition.
 \begin{theorem}\label{theorem5.2}
 		Let the vector $\boldsymbol{\mathcal{F}}=(\mathcal{F}_1,\cdots,\mathcal{F}_n)$, where $\mathcal{F}_i$ is one of the reinsurance families; proportional reinsurance family ($\mathcal{F}_P$), excess of loss reinsurance family ($\mathcal{F}_{XL}$) and limited excess of loss reinsurance family ($\mathcal{F}_{LXL}$).  Then, there exists an admissible reinsurance control $\boldsymbol{R}^* \in \boldsymbol{\mathcal{F}}$ such that $\pi({\cal{P}}^*, \boldsymbol{R}^*)$, the reinsurance band strategy  associated to ${\cal{P}}^*$ and $\boldsymbol{R}^*$,   is optimal. 
 \end{theorem}
\section{Numerical results}
For a  numerical solution of the value function and optimal control strategy, we use a method that is similar to the method described in Section 6.2 of \cite{azcue2014stochastic}. In fact, by use the finite difference method, we first solve Equation (\ref{valuefunction}) in the following way:  starting with  $$f_h(0)= a\,\, \text{and}\, f'_h(0)=\inf_{\boldsymbol{{R}} \in \boldsymbol{\mathcal{F}}}\frac{(\delta+(\sum_{i=1}^n\beta_i))a-(\sum_{i=1}^n\beta_i) P(\boldsymbol{{R}}(Y)=0)}{p_{\boldsymbol{{R}}}},\qquad a <\infty$$  and for $s=ih,\, i=1,2,\cdots$,  we approximate 
$\int_{0}^{x}f(x-\alpha) dF_R(\alpha)$
by  $$G_{\boldsymbol{R}}(s)=\sum_{\{j\le i\}}f_{h}((i-j)h)P\{(j-1)h<\boldsymbol{R}(Y)\le jh\}.$$
It is easy to show that $G_{\boldsymbol{R}}(s)$ converges to  $\int_{0}^{s}f(x-\alpha) dG_R(\alpha)$ as $h$ tends zero.
Then we define 
$f'_h(s)$ by
\begin{align}\label{3.24}
f'_{h}(s)&=\inf_{{\boldsymbol{R}}}\frac{(\delta+(\sum_{i=1}^n\beta_i)) h(f_{h}(s-h)-(\sum_{i=1}^n\beta_i) G_{{\boldsymbol{R}}}(s))}{{p_{{\boldsymbol{R}}}}},
\end{align}
and  set $f_h(s)=f_h(s-h)+hf'_h(s)$.
Now, if we set
\begin{align}
V_h(s)=
\begin{cases}	\frac{f_h(s)}{{f'_h(a_1^h)}}	& if \,\, s\le a_1^h,\\
s-a_1^h	+ \frac{f_h(a_1^h)}{{f'_h(a_1^h)}}& if \,\, s>a_1^h,
\end{cases}
\end{align}
where
$a_1^h={argmin}_s f'_h(s)$, 
then, 
for
$x\in [0,a_1^h]$ 
we have
$\sup_{\boldsymbol{{R}} \in \boldsymbol{\mathcal{F}}}\,\,{\cal{L}}_{\boldsymbol{R}}({V_h})(x)=0$ 
and
$\sup_{\boldsymbol{{R}} \in \boldsymbol{\mathcal{F}}}\,\,\Lambda_{\boldsymbol{R}}({V_h})(a_1^h)=0$.
 If 
$V_h$
is satisfied in (\ref{hjbeqV}) almost everywhere, we claim that the optimal strategy is as follows,
${\cal{A}}^h=\{a_1^h\}$, ${\cal{B}}^h=(a_1^h,\infty)$ 
and
${\cal{C}}^h=[0,a_1^h)$.
If $V_h$ does not satisfy in (\ref{hjbeqV}) almost everywhere, we consider the following function
\begin{align}
f_h^{(2)}(s,a_1^h,b_1,a_2)=
\begin{cases}	\frac{f_h(s)}{{f'_h(a_1^h)}}	& if \,\, s\le a_1^h,\\
s-a_1^h	+ \frac{f_h(a_1^h)}{{f'_h(a_1^h)}}	& if \,\, a_1^h<s\le b_1,\\
f_h^{(2)}(s)	& if \,\, b_1<s\le a_2,\\
f_h^{(2)}(a_2)-a_2+s & if \,\, a_2<s< \infty,
\end{cases}
\end{align}
where, for $s\in (b_1,a_2)$, 
$f_h^{(2)}(s)$
 is solved by numerical solution the (\ref{LV}), using the finite difference method and boundary condition
$\sup_{\boldsymbol{{R}} \in \boldsymbol{\mathcal{F}}}\,\,\Lambda_{\boldsymbol{R}}({f_h^{(2)}})(a_2)=0$.
Now, set
$V_h(s)= f_h^{(2)}(s,a_1^h,b_1^h,a_2^h)$
where
$$(b_1^h,a_2^h)=\underset{(b_1,a_2)}{argmax}  \{f_h^{(2)}(a_2)-a_2\}.$$
If
$V_h$
is satisfied in (\ref{hjbeqV}) almost everywhere, we claim that the optimal strategy is as follows,
${\cal{A}}^h=\{a_1^h,a_2^h\}$, ${\cal{B}}^h=(a_1^h,b_1^h]\bigcup (a_2^h,\infty)$
and
${\cal{C}}^h=[0,a_1^h)\bigcup (b_1,a_2^h) $. 
If $V_h$ does not satisfy in (\ref{hjbeqV}) almost everywhere,
we continue the process in the way described above.   The reinsurance and control strategy  obtained using the above algorithm are exhibited, respectively, by 
$\boldsymbol{R}^h$
and ${\cal{P}}^h=({\cal{A}}^h,{\cal{B}}^h,{\cal{C}}^h)$.
%%%%%%%%%%%%%%%%%
%%%%%%%%%%%%%%%%%%%%5
%%%%%%%%%%%%%%%%%%%%%%%%

\begin{theorem}
	If $V_h(x)\le{x}+{c}$, then the sequence $V_{h}$ converges to the unique viscosity $V$.
	\label{Theorem4.1}
\end{theorem}
\begin{figure}
	\centering
	\begin{subfigure}{.4\textwidth}
		\caption{}
		\centering
		\includegraphics[width=.8\linewidth]{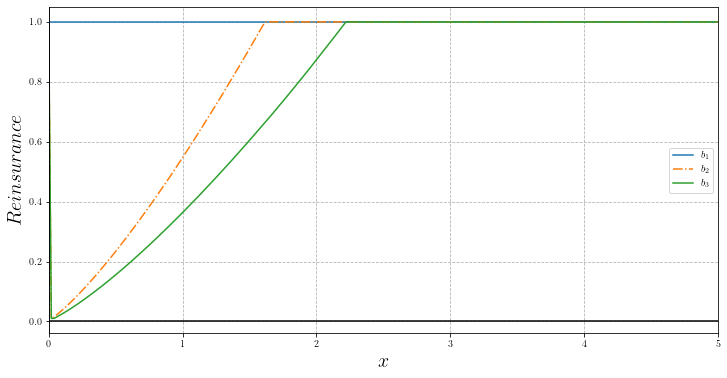}
		\label{one_p_RE}
	\end{subfigure}
	\begin{subfigure}{.4\textwidth}
		\caption{}
		\centering
		\includegraphics[width=.8\linewidth]{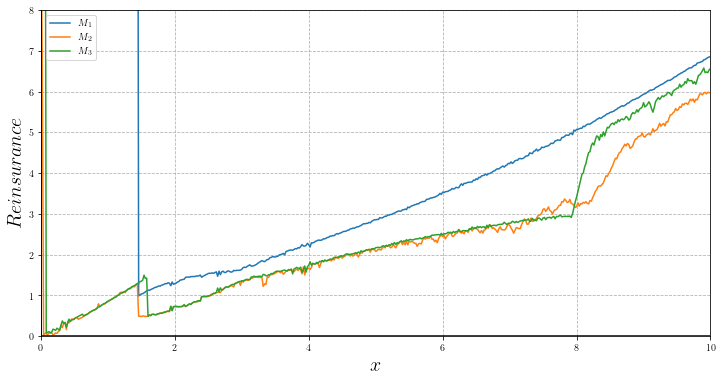}
		\label{multipleProRE}
	\end{subfigure}\\
	\begin{subfigure}{.4\textwidth}
		\caption{}
		\centering
		\includegraphics[width=.8\linewidth]{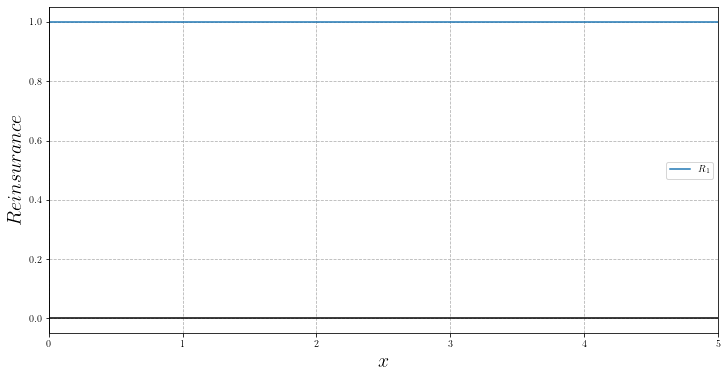}
		\label{oneXLRE}
	\end{subfigure}
	\begin{subfigure}{.4\textwidth}
		\caption{}
		\centering
		\includegraphics[width=.8\linewidth]{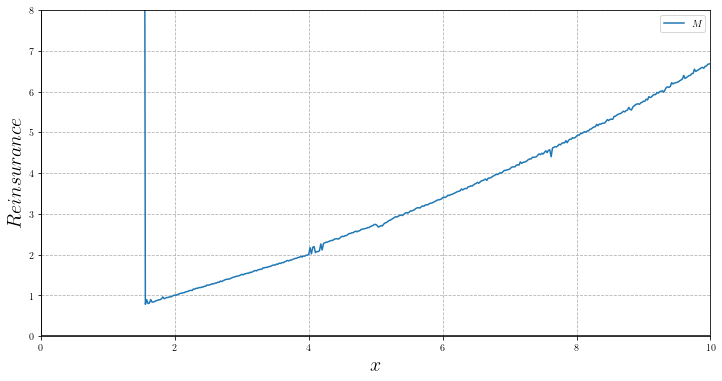}
		\label{multipleXLRE}
	\end{subfigure}
	\caption{ The numerical solution of the optimal reinsurances  with $h=0.02$, (a) The optimal results when three proportional reinsurances are used for three lines, (b) The optimal results when three XL reinsurances are used for three lines,  (c) The optimal result when one proportional reinsurances is used for three lines,  (d) The optimal results when one XL reinsurance is used for three lines }
	\label{optimal_R_indepen}
\end{figure}

Now, we obtain numerically some examples by using the above algorithm. 
\begin{example}
	Let insurance company has three lines of business such that it's  risk process has the  Thinning-dependence structure, defined in Remark \ref{remaek_example}; $F_i(x)=1-e^{-\lambda_i x}$, $i=1,2,3$, and $\lambda_1=0.5,\,\,\lambda_2=3,\,\, \lambda_3=2, \,\, \beta_1=8, \,\, \beta_2=4,\,\, \beta_3=5, \,\, \eta=3,  \,\,\eta_1=3.5$,
$	p_{11} = 1,\,\,
	p_{12} = 0.06,\,\,
	p_{13} = 0.05,\,\,	
	p_{21} = 0.03,\,\,
	p_{22} = 1,\,\,
	p_{23} = 0.01,\,\,
	p_{31} = 0.007,\,\,
	p_{32} = 0.005,\,\,
	p_{33} = 1.0$ and $\delta=0.3$. The reinsurance strategy in $i$th line is depicted by $R_i$. As was mentioned before, if $R_i \in \mathcal{F}_p$ then $R_i(y)=b_i(.) y$ and if $R_i \in \mathcal{F}_{XL}$ then $R_i(y)=\min(y,M_i(.))$, where $b_i(.)$ and $M_i(.)$ are functions of the company's capital.
	If the insurance company considers a reinsurance contract for three lines, the optimization issue will be equal with the uni-dimensional model scrutinized by \cite{azcue2014stochastic}. 
	Using the recently explained numerical method, the following results are gleaned,
	\begin{itemize}
		\item[(i)] if $R_i \in \mathcal{F}_p$, $i=1,\,2,\,3,$ then, ${\cal{P}}=\big(\{12.26\},\, (12.26, \infty),\,[0, 12.26)\big)$ ,
		\item[(ii)] if $R_i \in \mathcal{F}_{XL}$, $i=1,\,2,\,3,$ then, ${\cal{P}}=\big(\{10.44\},\, (10.44, \infty),\,[0, 10.44)\big)$,
		\item[(iii)] if $R_1=R_2=R_3 \in \mathcal{F}_p$,  then, ${\cal{P}}=\big(\{12.3\},\, (12.3, \infty),\,[0, 12.3)\big)$,
		\item[(iv)] if $R_1=R_2=R_3 \in \mathcal{F}_{XL}$,  then, ${\cal{P}}=\big(\{10.64\},\, (10.64, \infty),\,[0, 10.64)\big)$.
	\end{itemize}
\begin{figure}
	\centering
	\includegraphics[width=.8\linewidth]{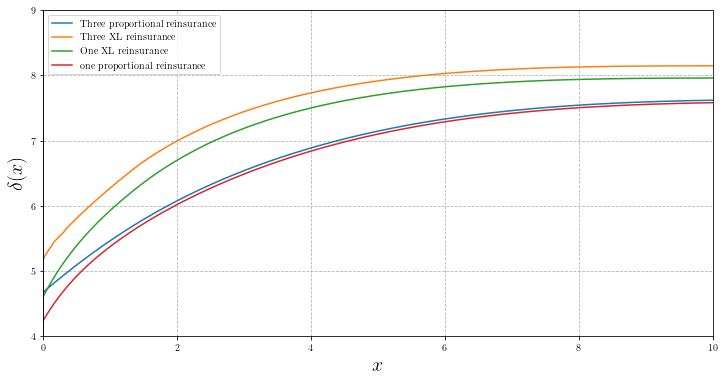}
	\caption{Survival functions}
	\label{multiplevaluefunction}
\end{figure}
\begin{figure}
	\centering
	\begin{subfigure}{.4\textwidth}
		\caption{}
		\centering
		\includegraphics[width=.8\linewidth]{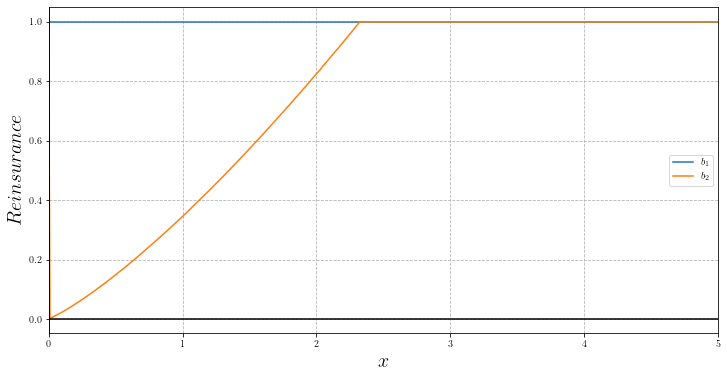}
		\label{one_p_RE}
	\end{subfigure}
	\begin{subfigure}{.4\textwidth}
		\caption{}
		\centering
		\includegraphics[width=.8\linewidth]{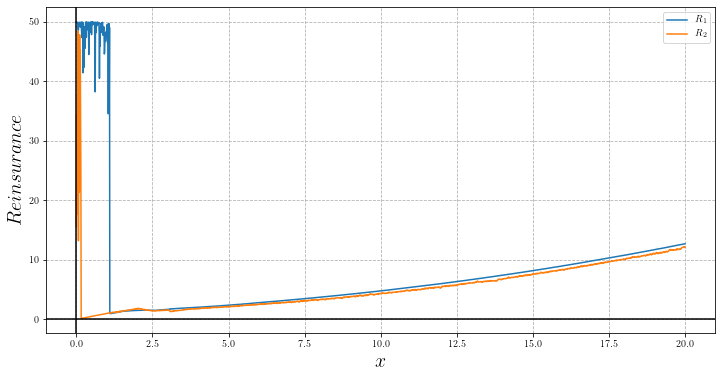}
		\label{multipleProRE}
	\end{subfigure}\\
	\begin{subfigure}{.4\textwidth}
		\caption{}
		\centering
		\includegraphics[width=.8\linewidth]{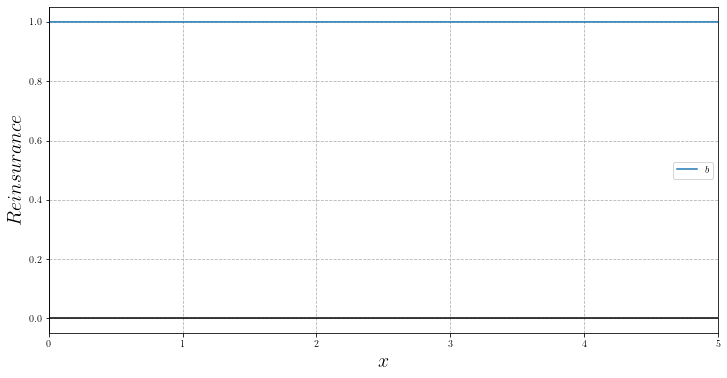}
		\label{oneXLRE}
	\end{subfigure}
	\begin{subfigure}{.4\textwidth}
		\caption{}
		\centering
		\includegraphics[width=.8\linewidth]{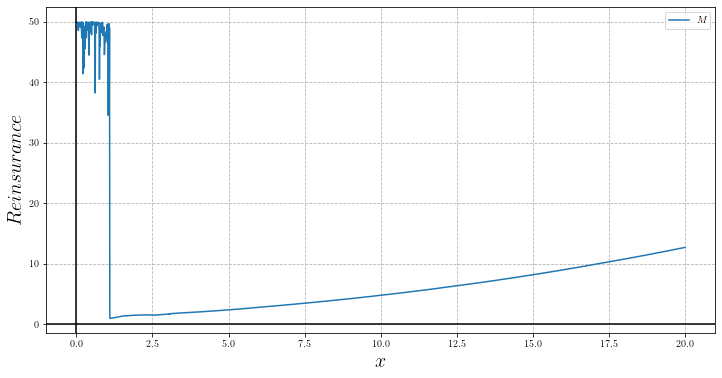}
		\label{multipleXLRE}
	\end{subfigure}
	\caption{ The numerical solution of the optimal reinsurances  with $h=0.01$ (a) The optimal results when two proportional reinsurances are used for two lines (b) The optimal results when two XL reinsurances are used for two lines  (c) The optimal result when one proportional reinsurances is used for two lines  (d) The optimal results when one XL reinsurance is used for two lines }
	\label{optimal_R_common}
\end{figure}
\begin{figure}
	\centering
	\includegraphics[width=.8\linewidth]{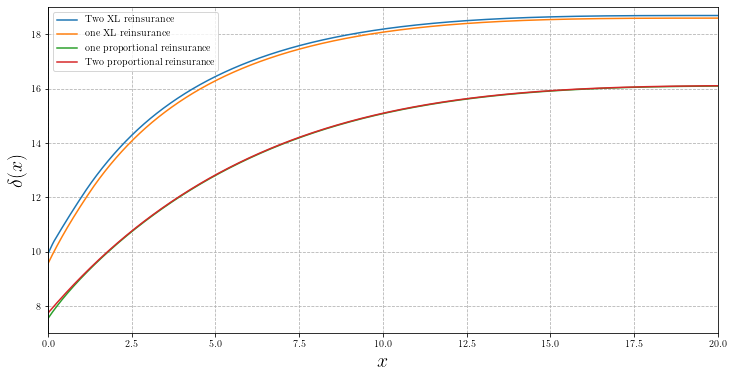}
	\caption{Survival functions}
	\label{multiplevaluefunction_commonshock}
\end{figure}
Also, optimization results for the value functions and reinsurance strategy are reported in Figures \ref{optimal_R_indepen} and \ref{multiplevaluefunction}.
\end{example}

\begin{example} (Common shock model)
	If $n = 2, m =3$, $p12 = p21 = 0$, $p31 = p32 = 1$ and $p11 = p22 = 1$, then $X ^{\pi}(t)$ defined in  (\ref{2.3'}) is
	the risk process with common shock for two dependent lines of business;
	$$	G^{\boldsymbol{{R}}}(\alpha)=\frac{\beta_1}{\sum_{i=1}^{3} \beta_i}P(R_1(U_1)\le \alpha)+\frac{\beta_2}{\sum_{i=1}^{3} \beta_i}P(R_2(U_2)\le \alpha)+\frac{\beta_3}{\sum_{i=1}^{3} \beta_i}P(R_1(U_1)+R_2(U_2)\le \alpha).$$
The following results are obtained by the numerical approach
	\begin{itemize}
		\item[(i)] if $R_i \in \mathcal{F}_p$, $i=1,\,2,$ then, ${\cal{P}}=\big(\{22.40\},\, (22.40, \infty),\,[0, 22.40)\big)$,
		\item[(ii)] if $R_i \in \mathcal{F}_{XL}$, $i=1,\,2,$ then, ${\cal{P}}=\big(\{19.83\},\, (19.83, \infty),\,[0, 19.83)\big)$,
		\item[(iii)] if $R_1=R_2=R_3 \in \mathcal{F}_p$,  then, ${\cal{P}}=\big(\{22.41\},\, (22.41, \infty),\,[0, 22.41)\big)$ and
		\item[(v)] if $R_1=R_2=R_3 \in \mathcal{F}_{XL}$,  then, ${\cal{P}}=\big(\{19.92\},\, (19.92, \infty),\,[0, 19.92)\big)$.
	\end{itemize}
Also, optimization results for the value functions and reinsurance strategy are reported in Figures \ref{optimal_R_common} and \ref{multiplevaluefunction_commonshock}.

\end{example}
\section{Conclusion}
%بررسی نحوه‌ی انتقال بخشی از ریسک شرکت بیمه به بیمه‌گر اتکایی و توزیع سود سهام بین سهام داران یک مساله‌ی است که در سال‌های اخیر مورد توجه محققین قرار گرفته است.  در حقیقت مساله این است که استراتژی اتکایی چگونه باشد تا مقدارتنزیلی سود سهام توزیع شده بین سهام‌داران (تابع هدف) تا زمان ورشکستگی ماکزیمم شود؟
To solve this optimization issue, an HJB equation associated with the value function $V(.)$ defined in (\ref{valuefunction}) is adapted. 
Usually, the value  function does not have the smoothness properties required for interpreting it as a solution for corresponding HJB equation in the classical sense, but it is satisfying in this equation in a weaker concept.  The optimal survival function is characterized  as the smallest super-viscosity solution of the HJB equation.
Unfortunately, obtaining a closed form for the value function or the control strategy in the issue discussed in this paper is  complicated or impossible. Therefore, it was more practical to adopt a numerical solution.  For constructing a numerical solution, the FDM has been employed because the convergence of a numerical solution to the value function can be proved through the techniques prevalent in the literature. The convergent findings are displayed in section 4. 
The results  of the present paper give the insurance companies this opportunity to share their risk with the reinsurers. In section 5,  examples reveal that using this approach, the survival function will be increased. To sum up, with the implementation of this dynamic method for drawing the vector of the reinsurance contracts, the value function $V(.)$ might augment significantly.

\section*{Appendix A}

  %%%%%%%%%%%%%%%%%%%%%%%%%%%%%%%%%%%%%%
  {\bf  Proof of Proposition  \ref{prop5.5}}\,\, Let us define
  \begin{align*}
  {\cal{T}}(f)(x)=E_x\bigg[\int_{0}^{\tau_1}e^{-\delta t}dD_s+e^{-\delta \tau_1}f\big(X_{\tau_1}^{\pi}\big)\bigg],
  \end{align*}
  where $\tau_1$ is the time of the first claim and $\pi=\pi({\cal{P}},\boldsymbol{r})=\{(\boldsymbol{R}^x, D^x)\in \Pi_{x}^{\pi}, x\ge 0\}$ defined in Definition \ref{defi5.3}. It is obvious that  ${\cal{T}}(V_{\pi})=V_{\pi}$, i.e., $V_{\pi}$ is a fixed point of ${\cal{T}}$. Let any $M\in {\cal{A}}\bigcup {\cal{B}}$, define the following complete metric space
  $${\cal{B}}_M=\{f:[0,M]\longrightarrow [0,\infty) \,\,\text{Borel-measurable and bounded}\}$$
  with the metric $d(f_1,f_2)=\sup_{x\in[0,M]}|f_1(x)-f_2(x)|$. For any $x\le M$, we have that $X_t^{\pi}\le M$ and so ${\cal{T}}$ is well defined and bounded in $[0,M]$. It is easy to see that 
  $$|{\cal{T}}(f_1)(x)-{\cal{T}}(f_2)(x)|\le \frac{\sum_{i=1}^{m}\beta_i}{\delta+\sum_{i=1}^{m}\beta_i}d(f_1,f_2).$$
  Therefore ${\cal{T}}:{\cal{B}}_M\longrightarrow {\cal{B}}_M$ is a contraction with modulus $\frac{\sum_{i=1}^{m}\beta_i}{\delta+\sum_{i=1}^{m}\beta_i}<1$ and so, by the contraction mapping theorem, ${\cal{T}}$ has a unique fixed point.
  
  Now, for complete proof, it is enough to show that ${\cal{T}}(W)=W$. If $x\in {\cal{A}}^*$, then 
  $${\cal{T}}(x)=\frac{1}{c+\sum_{i=1}^{m}\beta_i}\big(p_{\boldsymbol{R}}^*-\sum_{i=1}^{m}\beta_i\int_{0}^{x}V(x-\alpha)dG^{\boldsymbol{R}}(\alpha)\big).$$
  Since $\Lambda(x)=0$ we have ${\cal{T}}(x)=\Lambda(x)$. If $x \in {\cal{B}}$, we have  $x_0=\max\{y<x\,\, \text{and}\,\,y\notin {\cal{B}}  \}\in {\cal{A}}$   and so 
  $${\cal{T}}(W)(x)=x-x_0+W(x_0)=W(x_0).$$
  For $x\in {\cal{C}}$, consider $x_1=\min\{y>x\,\,\text{and}\,\,y\notin {\cal{C}} \}\in {\cal{A}}$ and 
  $$X(t)=x+\int_{0}^{t}p_{\boldsymbol{R}_s}^*ds.$$
  We can find $t_1>0$ such that $X(t_1)=x_1$ and for $t\in [0,t_1),$ $X(t)\in {\cal{C}}$. So, ${\cal{T}}(W)(x_1)=W(x_1)$ and 
  \begin{align*}
  {\cal{T}}(x)&=e^{-(c+\sum_{i=1}^{m}\beta_i)t_1}W(x_1)+\int_{0}^{t_1}\bigg(\int_{0}^{x(t)}W(x(t)-\alpha)dG^{\boldsymbol{{R}}}(\alpha)\bigg)(\sum_{i=1}^{m}\beta_i)e^{-(c+\sum_{i=1}^{m}\beta_i)}dt\\
  &=e^{-(c+\sum_{i=1}^{m}\beta_i)t_1}W(x_1)+\int_{0}^{t_1}\big(p_{\boldsymbol{{R}}}W'(x(t))-(c+\sum_{i=1}^{m}\beta_i)W(X(t))\big)e^{-(c+\sum_{i=1}^{m}\beta_i)}dt\\
  &=e^{-(c+\sum_{i=1}^{m}\beta_i)t_1}W(x_1)+\int_{0}^{t_1}(W(x(t))e^{-(c+\sum_{i=1}^{m}\beta_i)})'dt\\
  &=W(x).
  \end{align*}
  $\blacksquare$
  %%%%%%%%%%%%%%%%%%%%%%%%%%%
  
  Suppose 
  $$X(t)=x+p_R-\sum_{i=1}^{N_t}R(Z_i)-D_t.$$
  The HJB equation related to in $V(x)=\sup_{R \in {\mathcal{F}}}\int_{0}^{\tau^R}e^{-\delta s}dL_s$ is as follows
  
$$
  \max\{1-V'(x),{sup}\,\,{{\cal{L}}}_{R}(V)(x)\}=0,
$$
  where
  $$ {{\cal{L}}}_{\boldsymbol{R}}(V)(x)=p_{R}V'(x)-(c+\sum_{i=1}^m \beta_i)V(x)+(\sum_{i=1}^m \beta_i)\int_{0}^{x}V(x-\alpha)dG^{R}(\alpha),$$
  and $G^{R}$ is a distribution function under a reinsurance strategy. If that is so, our results are in harmony with the finding of \cite{azcue2005optimal} and \cite{azcue2014stochastic}. In the present paper, the HJB equation pertained to $V(.)$ is as follows
 $$
  \max\{1-V'(x),\underset{\boldsymbol{\mathcal{F}}}{sup}\,\,{{\cal{L}}}_{\boldsymbol{R}}(V)(x)\}=0,
$$
where
$$ {{\cal{L}}}_{\boldsymbol{R}}(V)(x)=p_{\boldsymbol{R}}V'(x)-(c+\sum_{i=1}^m \beta_i)V(x)+(\sum_{i=1}^m \beta_i)\int_{0}^{x}V(x-\alpha)dG^{\boldsymbol{R}}(\alpha),$$ 
  and 
  $G^{\boldsymbol{R}}(\alpha)$ is defined in (\ref{Gdistribution}). Therefore, the proof of Proposition 3.8, Proposition 4.2 and Proposition 5.1 from  \cite{azcue2005optimal} hold water for Theorem \ref{V-viscosity}, Proposition \ref{prop4.3.} and Proposition \ref{prop4.4}. Moreover, Proposition \ref{prop5.7} and Theorem \ref{theorem5.2} are akin to Proposition 5.7 and Theorem 5.2 by \cite{azcue2014stochastic}; notwithstanding, proving these theorems depends on the structure of reinsurance strategy which will be discussed accordingly.
  
  %%%%%%%%%%%%%%%%%%%%%%%%%%%%%
  Before going in the proof of Proposition \ref{prop5.7} and Theorem \ref{theorem5.2}, we need the following lemmas.
  \begin{lemma}\label{remark5.5}
  {For any} $x\le \tilde{x}$, define
  	$$\Pi_{x,\tilde{x}}=\{\pi\in \Pi_{x} \,\,\text{such that }\,\, X_t^{\pi}<\tilde{x}\,\,\text{for all}\,\, t \ge 0\}.$$
  	If $\tilde{x}$ satisfying either $sup_{\boldsymbol{{R}} \in \boldsymbol{\mathcal{F}}}\Lambda_{\boldsymbol{R}}(W)(\tilde{x})=0$ or $V'(\tilde{x})=1$, then 
  	$$V(x)=\underset{\boldsymbol{\pi}\in \Pi_x}{\sup}\,\,V^{\pi}(x).$$
  \end{lemma}
 \begin{lemma}\label{lemma3.2.'}
	Let the vector $\boldsymbol{\mathcal{F}}=(\mathcal{F}_1,\cdots,\mathcal{F}_n)$, where $\mathcal{F}_i$ is one of the reinsurance families $\mathcal{F}_P$, $\mathcal{F}_{XL}$ and $\mathcal{F}_{LXL}$. If u is positive and continuously differentiable,
then the function $sup_{\boldsymbol{{R}} \in \boldsymbol{\mathcal{F}}}\Lambda_{\boldsymbol{R}}(u)$ is right continuous and upper semicontinuous. Moreover,  $sup_{\boldsymbol{{R}} \in \boldsymbol{\mathcal{F}}}\Lambda_{\boldsymbol{R}}(u)\le 0.$ 
\end{lemma}
{\bf Proof}  For proving this lemma, we must prove that
$sup_{\boldsymbol{{R}} \in \boldsymbol{\mathcal{F}}}\Lambda_{\boldsymbol{R}}(u)$ is right upper semicontinuous, left lower semicontinuous and right lower semicontinuous. Right continuous is follows from right upper semicontinuous and right lower semicontinuous, and upper semicontinuous is follows from right upper semicontinuous and left upper semicontinuous.
Let us prove first that $sup_{\boldsymbol{{R}} \in \boldsymbol{\mathcal{F}}}\Lambda_{\boldsymbol{R}}(u)$  is left upper semicontinuous. For assumed $x_0$, $x_k \nearrow x_0$, consider reinsurance strategies ${\boldsymbol{R}}^{(k)}\in \boldsymbol{\mathcal{F}}$ such that
\begin{align}\label{A.5}
\sup_{\boldsymbol{R}\in \boldsymbol{\mathcal{F}}} \Lambda_{\boldsymbol{R}}(u)(x_k)\le \Lambda_{{\boldsymbol{R}}^{(k)}}(u)(x_k)+\frac{1}{k}.
\end{align}
Then, the following is straightforward,
\begin{align*}
\Lambda_{{\boldsymbol{R}}^{(k)}}(u)(x_0)=&\Lambda_{{\boldsymbol{R}}^{(k)}}(u)(x_k)-(c+\sum_{i=1}^m\beta_i)(u(x_0)-(x_k))+\sum_{i=1}^m\beta_i\int_{x_k}^{x_0}u(x_0-\alpha) dF_{{\boldsymbol{R}}}^{(k)}(\alpha)\\
&+\sum_{i=1}^m\beta_i\int_{0}^{x_k}(u(x_0-\alpha)-u(x_k-\alpha)) dF_{{\boldsymbol{R}}}^{(k)}(\alpha)\\
\ge &\Lambda_{{\boldsymbol{R}}^{(k)}}(u)(x_k)-(c+\sum_{i=1}^m\beta_i)(u(x_0)-u(x_k))+\sum_{i=1}^m\beta_i\int_{0}^{x_k}(u(x_0-\alpha)-u(x_k-\alpha)) dF_{{{\boldsymbol{R}}}^{(k)}}(\alpha).
\end{align*} 
So, we have,
\begin{align}\label{A.6}
\limsup_{k\rightarrow \infty}\Lambda_{{\boldsymbol{R}}^{(k)}}(u)(x_0)\ge \limsup_{k\rightarrow \infty}\Lambda_{{\boldsymbol{R}}^{(k)}}(u)(x_k).
\end{align}
Then, from (\ref{A.5}) and (\ref{A.6}), the following result can be derived; 
$$\sup_{\boldsymbol{R}\in \boldsymbol{\mathcal{F}}} \Lambda_{\boldsymbol{R}}(u)(x_0)\ge \limsup_{k\rightarrow \infty}\Lambda_{{\boldsymbol{R}}^{(k)}}(u)(x_0)\ge \limsup_{k\rightarrow \infty}\left(\sup_{\boldsymbol{R}\in \boldsymbol{\mathcal{F}}} \Lambda_{\boldsymbol{R}}(u)(x_k)\right).$$
Consequently, the following relation is dominant
$\limsup_{x\rightarrow x_0^-}\left(\sup_{\boldsymbol{R}\in \boldsymbol{\mathcal{F}}} \Lambda_{\boldsymbol{R}}(u)(x)\right)\le \sup_{\boldsymbol{R}\in \boldsymbol{\mathcal{F}}} \Lambda_{\boldsymbol{R}}(u)(x_0).$

Now, Let us prove first that $sup_{\boldsymbol{{R}} \in \boldsymbol{\mathcal{F}}}\Lambda_{\boldsymbol{R}}(u)$  is right upper semicontinuous. We must to show that  $\limsup_{x\rightarrow x_0^+}\left(\sup_{\boldsymbol{R}\in \boldsymbol{\mathcal{F}}} \Lambda_{\boldsymbol{R}}(u)(x)\right)\le \sup_{\boldsymbol{R}\in \boldsymbol{\mathcal{F}}} \Lambda_{\boldsymbol{R}}(u)(x_0).$
Given any sequence $x_k\searrow x_0$, take the reinsurance strategies ${\boldsymbol{R}}^{(k)}\in \boldsymbol{\mathcal{F}}$ such that 
$$
\sup_{\boldsymbol{R}\in \boldsymbol{\mathcal{F}}} \Lambda_{\boldsymbol{R}}(u)(x_k)\le \Lambda_{{\boldsymbol{R}}^{(k)}}(u)(x_k)+\frac{1}{k}.
$$
If one of the reinsurance contracts is LXL reinsurance, for example $\mathcal{F}_1=\mathcal{F}_{LXL}$, 
take $$\bar{{\boldsymbol{R}}}^{(k)}=(\bar{R}_1^{(k)},\bar{R}_2^{(k)},\cdots , \bar{R}_n^{(k)}) \in \boldsymbol{\mathcal{F}}$$ such that
\begin{align*}
\bar{R}_1^{(k)}(\alpha)=
\begin{cases}	R_1^{(k)}(\alpha)	& if \,\, R_1^{(k)}(\alpha)=\alpha\,\, for\,\, all\,\, \alpha \\
R_1^{(k)}(\alpha)	& if \,\, R_1^{(k)}(\alpha)=a_k\wedge \alpha+(\alpha-L-a_k)^+\,\, with\,\, a_k \notin (x_0,x_k) \\
\alpha \wedge x_0 +(\alpha-L-x_0)^+	& if \,\, R_1^{(k)}(\alpha)=a_k\wedge \alpha+(\alpha-L-a_k)^+\,\, with\,\, a_k \in (x_0,x_k) 
\end{cases}
\end{align*}
and $\bar{R}_2^{(k)}=R_2^{(k)}, \cdots, \bar{R}_n^{(k)}=R_n^{(k)}$.
If
$a_k \le x_0$ 
then
\begin{align*}
&\int_{0}^{\infty} (u(x_k-\alpha)-u(x_0-\alpha))dF_{\bar{{\boldsymbol{R}}}^{(k)}}(\alpha)\\
&= \int_{0}^{\infty} (u(x_k-\alpha)-u(x_0-\alpha))dF_{{\boldsymbol{R}}^{(k)}}(\alpha)\\
&\le \sup_{x \in [0,A]}|u'(x)|(y-x_0).
\end{align*}
Suppose $S^{n-1}=\{2, \cdots,n\}$, and $S_i^{n-1}=\{A_{ij}^{n-1}; j=1,\cdots,\tiny{\begin{pmatrix}n-1\\i\end{pmatrix}}  \}$, where $A_{ij}^{n-1}$ a subset of a set $S^{n-1}$, with exactly $i$ elements, then we have
\begin{align*}
G^{\boldsymbol{R}}(\alpha)=&\sum_{j=1}^{n}\sum_{k=1}^{\tiny{\begin{pmatrix}n\\j\end{pmatrix}}}\bigg(\sum_{i=1}^{m}\frac{\beta_i}{\sum_{i=1}^{m}\beta_i}\prod_{z\in A_{jk}^{n}}p_{iz}\prod_{z\in S^{n}-A_{jk}^{n}}(1-p_{iz})\bigg)F_{{\boldsymbol{R}}_{A_{jk}^n}}(\alpha)\\
=&\sum_{i=1}^{m}\frac{\beta_i}{\sum_{i=1}^{m}\beta_i}p_{i1}\prod_{z\neq 1}(1-p_{iz})F_{R_1}(\alpha_1)
\\
&+\sum_{j=1}^{n-1}\sum_{k=1}^{\tiny{\begin{pmatrix}n-1\\j\end{pmatrix}}}\bigg(\sum_{i=1}^{m}\frac{\beta_i}{\sum_{i=1}^{m}\beta_i}p_{i1}\prod_{z\in A_{jk}^{n-1}}p_{iz}\prod_{z\in S^{n-1}-A_{jk}^{n-1}}(1-p_{iz})\bigg)\int_{0}^{\infty}F_{{\boldsymbol{R}}_{A_{jk}^{n-1}}}(\alpha-R_1(\alpha_1))dF_{R_1}(\alpha_1)\\
&+\sum_{j=1}^{n-1}\sum_{k=1}^{\tiny{\begin{pmatrix}n-1\\j\end{pmatrix}}}\bigg(\sum_{i=1}^{m}\frac{\beta_i}{\sum_{i=1}^{m}\beta_i}(1-p_{i1})\prod_{z\in A_{jk}^{n-1}}p_{iz}\prod_{z\in S^{n-1}-A_{jk}^{n-1}}(1-p_{iz})\bigg)F_{{\boldsymbol{R}}_{A_{jk}^{n-1}}}(\alpha).\end{align*}
Let us define $ {{\boldsymbol{R}}}_{n-1}^{(k)}=(R_2^{(k)},\cdots, R_n^{(k)})$. If
$a_k \in (x_0,x_k)$, 
then
\begin{align*}
&\int_{0}^{\infty}  u(x_k-\alpha)dG^{{{\boldsymbol{R}}}^{(k)}}(\alpha)- \int_{0}^{\infty} u(x_0-\alpha)dG^{\bar{{\boldsymbol{R}}}^{(k)}}(\alpha)\\
=&
\sum_{i=1}^{m}\frac{\beta_i}{\sum_{i=1}^{m}\beta_i}p_{i1}\prod_{z\neq 1}(1-p_{iz})\times
\bigg( \int_{0}^{x_0} \big( u(x_k-R_1^{(k)}(\alpha))-u(x_0-R_1^{(k)}(\alpha))dF_1(\alpha)\\
&+\int_{x_0}^{a_k} \big( u(x_k-\alpha-(\alpha-L-a_k)^+)-u(-(\alpha-L-x_0)^+)\big)dF_1(\alpha)\\
&+ \int_{a_k}^{x_k} \big( u(x_k-a_k-(\alpha-L-a_k)^+)-u(-(\alpha-L-x_0)^+)\big)dF_1(\alpha)
\\
&+ \int_{x_k}^{\infty} \big( u(x_k-\alpha_k-(\alpha-L-a_k)^+)-u(-(\alpha-L-x_0)^+)\big)dF_1(\alpha)\bigg)
\\
&+
\sum_{j=1}^{n-1}\sum_{k=1}^{\tiny{\begin{pmatrix}n-1\\j\end{pmatrix}}}\sum_{i=1}^{m}\frac{\beta_i}{\sum_{i=1}^{m}\beta_i}p_{i1}\prod_{z\in A_{jk}^{n-1}}p_{iz}\prod_{z\in S^{n-1}-A_{jk}^{n-1}}(1-p_{iz})\\
&\times\bigg( \int_{0}^{x_0} \int_{0}^{\infty}\big( V(x_k-R_1^{(k)}(\alpha_1)-\alpha_2)-u(x_0-R_1^{(k)}(\alpha_1)-\alpha_2)dF_{{\boldsymbol{R}}^{(k)}_{A_{jk}^{n-1}}}(\alpha_2)dF_1(\alpha_1)\\
&+ \int_{x_0}^{a_k} \int_{0}^{\infty}\big( u(x_k-\alpha_1-(\alpha_1-L-a_k)^+-\alpha_2)-u(-(\alpha_1-L-x_0)^+-\alpha_2)dF_{{\boldsymbol{R}}^{(k)}_{A_{jk}^{n-1}}}(\alpha_2)dF_1(\alpha_1)\\
&+ \int_{a_k}^{x_k} \int_{0}^{\infty}\big( u(x_k-\alpha_k-(\alpha_1-L-a_k)^+-\alpha_2)-u(-(\alpha_1-L-x_0)^+-\alpha_2)dF_{{\boldsymbol{R}}^{(k)}_{A_{jk}^{n-1}}}(\alpha_2)dF_1(\alpha_1)\\
&+ \int_{x_k}^{\infty} \int_{0}^{\infty}\big( u(x_k-\alpha_k-(\alpha_1-L-a_k)^+-\alpha_2)-u(-(\alpha_1-L-x_0)^+-\alpha_2)dF_{{\boldsymbol{R}}^{(k)}_{A_{jk}^{n-1}}}(\alpha_2)dF_1(\alpha_1)\bigg)\\
&\le 2 \sup_{x \in [0,x_k]}|u'(x)|(x_k-x_0)\\
&+\sup_{\alpha \in [x_0,a_k]}(u(x_k-\alpha-(\alpha-L-a_k)^+)-u(-(\alpha-L-x_0)^+))p(x_0\le Z_1 \le a_k)\\
&+\sup_{\alpha \in [a_k,x_k]}(u(x_k-a_k-(\alpha-L-a_k)^+)-u(-(\alpha-L-x_0)^+))p(a_k\le  Z_1\le x_k)\\
&+\sup_{\alpha \in [x_k,\infty)}(u(x_k-a_k-(\alpha-L-a_k)^+)-u(-(\alpha-L-x_0)^+))p(  Z_1 \ge x_k)\\
&+\sup_{\alpha \in [x_0,a_k]}\sup_{\alpha_2 \in [0,x_k]}(u(x_k-\alpha-(\alpha-L-a_k)^+-\alpha_2)-u(-(\alpha-L-x_0)^+-\alpha_2))p(x_0\le Z_1 \le a_k)\\
&+\sup_{\alpha \in [a_k,x_k]}\sup_{\alpha_2 \in [0,x_k]}(u(x_k-a_k-(\alpha-L-a_k)^+-\alpha_2)-u(-(\alpha-L-x_0)^+)-\alpha_2)p(a_k\le  Z_1\le x_k)\\
&+\sup_{\alpha \in [x_k,\infty)}\sup_{\alpha_2 \in [0,x_k]}(u(x_k-a_k-(\alpha-L-a_k)^+-\alpha_2)-u(-(\alpha-L-x_0)^+-\alpha_2))p(  Z_1 \ge x_k)
\end{align*}
According to the property of  right-continuously of  distribution function; if 
$x_k\searrow x_0$,
then
$$p(x_0\le Z_1 \le a_k)\longrightarrow 0\,\,\,\,
\text{and}\,\,\,\,
p(a_k\le Z_1 \le x_k)\longrightarrow 0.$$
Now, the term 
$\sup_{\alpha \in [x_k,\infty)}(u(x_k-a_k-(\alpha-L-a_k)^+)-u(-(\alpha-L-x_0)^+))$ 
should become the focus of attention. In this case, there are two situations as outlined below:
\begin{itemize}
	\item [(I)]If there is a finite value $ m$ satisfying the following,
	\begin{align*}
	&\sup_{\alpha \in [x_k,\infty)}(u(x_k-a_k-(\alpha-L-a_k)^+)-u(-(\alpha-L-x_0)^+))\\
	&=u(x_k-a_k-(m-L-a_k)^+)-u(-(m-L-x_0)^+).
	\end{align*}
	So, 
	\begin{align*}
	& \lim_{x_k\searrow x_0}\sup_{\alpha \in [x_k,\infty)}(u(x_k-a_k-(\alpha-L-a_k)^+)-u(-(\alpha-L-x_0)^+))\\
	&=u(-(m-L-x_0)^+)-u(-(m-L-x_0)^+)=0.
	\end{align*}  
	\item[(II)]  
	If
	\begin{align*}
	& \sup_{\alpha \in [x_k,\infty)}(u(x_k-a_k-(\alpha-L-a_k)^+)-u(-(\alpha-L-x_0)^+))\\
	&=\lim_{\alpha \rightarrow \infty}u(x_k-a_k-(\alpha-L-a_k)^+)-u(-(\alpha-L-x_0)^+)
	\end{align*}
	then
	$$
	\sup_{\alpha \in [x_k,\infty)}(u(x_k-a_k-(\alpha-L-a_k)^+)-u(-(\alpha-L-x_0)^+))
	=u(x_k-a_k)-u(0).
	$$
	So,
	\begin{align*}
	&\lim_{x_k\searrow x_0}\sup_{\alpha \in [x_k,\infty)}(u(x_k-a_k-(\alpha-L-a_k)^+)-u(-(\alpha-L-x_0)^+))\\
	&=\lim_{x_k\searrow x_0}(u(x_k-a_k)-u(0))=0.
	\end{align*}  
\end{itemize}
With a similar
argument we can obtain
$$
\lim_{x_k\searrow x_0}\sup_{\alpha \in [x_k,\infty)}\sup_{\alpha_2 \in [0,x_k]}(u(x_k-a_k-(\alpha-L-a_k)^+-\alpha_2)-u(-(\alpha-L-x_0)^+-\alpha_2))=0.
$$
It should be noted that  
$R_1^{(k)}(\alpha)-(x_k-x_0)\le \bar{R}_1^{(k)}(\alpha)\le R_1^{(k)}(\alpha).$
So 
$p_{\bar{R}_1^{(k)}} \nearrow p_{R_1^{(k)}}$. 
Thus
\begin{align*}
&sup_{\boldsymbol{{R}} \in \boldsymbol{\mathcal{F}}}\Lambda_{\boldsymbol{R}}(u)(x_k)- sup_{\boldsymbol{{R}} \in \boldsymbol{\mathcal{F}}}\Lambda_{\boldsymbol{R}}(u)(x_0)\le \Lambda_{\boldsymbol{R}^{(k)}}(u)(x_k)-\Lambda_{{\boldsymbol{R}}^{(k)}}(v)(x_0)+\epsilon\\
&\le 2(c+\sum_{i=1}^m\beta_i) \sup_{x\in [0,A]}|u'(x)|(x_k-x_0)+2 \sup_{x \in [0,A]}|u'(x)|(x_k-x_0)\\
&+\sup_{\alpha \in [x_0,a_k]}(u(x_k-\alpha-(\alpha-L-a_k)^+)-u(-(\alpha-L-x_0)^+))p(x_0\le Z_1 \le a_k)\\
&+\sup_{\alpha \in [a_k,x_k]}(u(x_k-a_k-(\alpha-L-a_k)^+)-u(-(\alpha-L-x_0)^+))p(a_k\le  Z_1\le x_k)\\
&+\sup_{\alpha \in [x_k,\infty)}(u(x_k-a_k-(\alpha-L-a_k)^+)-u(-(\alpha-L-x_0)^+))p(  Z_1 \ge x_k)\\
&+\sup_{\alpha \in [x_0,a_k]}\sup_{\alpha_2 \in [0,x_k]}(u(x_k-\alpha-(\alpha-L-a_k)^+-\alpha_2)-u(-(\alpha-L-x_0)^+-\alpha_2))p(x_0\le Z_1 \le a_k)\\
&+\sup_{\alpha \in [a_k,x_k]}\sup_{\alpha_2 \in [0,x_k]}(u(x_k-a_k-(\alpha-L-a_k)^+-\alpha_2)-u(-(\alpha-L-x_0)^+)-\alpha_2)p(a_k\le  Z_1\le x_k)\\
&+\sup_{\alpha \in [x_k,\infty)}\sup_{\alpha_2 \in [0,x_k]}(u(x_k-a_k-(\alpha-L-a_k)^+-\alpha_2)-u(-(\alpha-L-x_0)^+-\alpha_2))p(  Z_1 \ge x_k)
\end{align*}
and so we get that $sup_{\boldsymbol{{R}} \in \boldsymbol{\mathcal{F}}}\Lambda_{\boldsymbol{R}}(u)(.)$ is right upper semicontinuous. 
The proof for the case $\mathcal{F}_1=\mathcal{F}_{XL}$ and $\mathcal{F}_1=\mathcal{F}_{p}$ are simpler,  we therefore omit them.
Now, repeating the arguments presented in the proof of  Proposition 7.4 of \cite{azcue2005optimal} (replacing   $\boldsymbol{R}=(\boldsymbol{R}_t)_{t \ge 0}=({R_1}_t,\cdots, {R_n}_t)_{t \ge 0}$ with  $\bar{R}$), the right lower semicontinuous is obtained.
$\square$
\begin{lemma}
	\label{lemma5.6}
	For any $\boldsymbol{{R}}\in \boldsymbol{\mathcal{F}}$,  let us define the function
	$$ \tilde{{\cal{L}}}_{\boldsymbol{R}}(V)(x)=p_{\boldsymbol{R}}\hat{V}(x)-(c+\sum_{i=1}^m \beta_i)V(x)+(\sum_{i=1}^m \beta_i)\int_{0}^{x}V(x-\alpha)dG^{\boldsymbol{R}}(\alpha)$$
	where
	$$\hat{V}(x)=\inf_{\boldsymbol{{R}}}\frac{(c+\sum_{i=1}^m \beta_i)V(x)-(\sum_{i=1}^m \beta_i)\int_{0}^{x}V(x-\alpha)dG^{\boldsymbol{R}}(\alpha)}{p_{\boldsymbol{R}}}.$$
Then, we have the following result;
\begin{itemize}
	\item[(a)] 
	$\hat{V}$ is well defined and Borel measurable, $\hat{V}\ge 1$, and 	$$sup_{\boldsymbol{{R}} \in \boldsymbol{\mathcal{F}}}\,\,{\tilde{{\cal{L}}}}_{\boldsymbol{R}}({V})(x)=0.$$ 
	\item[(b)] If $x\in {\cal{C}}^*$:  $\hat{V}(x)>1$, and $\hat{V}(x)=V'(x)$ where $V$ is differentiable (i.e., a.e. in $x\in {\cal{C}}$).
	\item[(c)] If $\hat{V}(x)=1$ then $x\in {\cal{A}}^*$.
\end{itemize}
\end{lemma}
\begin{lemma}\label{prop5.6}
	Let the vector $\boldsymbol{\mathcal{F}}=(\mathcal{F}_1,\cdots,\mathcal{F}_n)$, where $\mathcal{F}_i$ is one of the reinsurance families $\mathcal{F}_P$, $\mathcal{F}_{XL}$ and $\mathcal{F}_{LXL}$.  Then, there exists a $\boldsymbol{R}_x^* \in \boldsymbol{\mathcal{F}}$ such that
	$$sup_{\boldsymbol{{R}} \in \boldsymbol{\mathcal{F}}}\,\,{\tilde{{\cal{L}}}}_{\boldsymbol{R}}({V})(x)={\tilde{{\cal{L}}}}_{\boldsymbol{R}_x^*}({V})(x).$$ 
\end{lemma}
{\bf  Proof }\,\,It is enough to show that, there exists  a $\boldsymbol{R}^*=({R}_1^*,\cdots,{R}_n^*)\in \boldsymbol{\mathcal{F}}$, where ${R}_i^*$ is one of the reinsurance families ${R}_P$, ${R}_{XL}$ and ${R}_{LXL}$, such that the maximum of 
\begin{align*}
g(V,x,\boldsymbol{R})&=(1+\eta_1)\int_{0}^{\infty}\alpha dG^{\boldsymbol{R}}(\alpha)+\int_{0}^{x}V(x-\alpha)dG^{\boldsymbol{R}}(\alpha)\\
&=\int_{0}^{\infty}\big((1+\eta_1)\alpha +V(x-\alpha)\big)dG^{\boldsymbol{R}}(\alpha)
\end{align*}
is attained at $\boldsymbol{R}^*$. Let us assume
\begin{align*}
R_i(\alpha)&= b_i\alpha, \,\,i=1,\cdots, k_1\\
R_i(\alpha)&= a_i\wedge\alpha, \,\,i=k_1+1,\cdots, k_1+k_2\\
R_i(\alpha)&= a_i\wedge\alpha+(\alpha-a_i-L_i)^+,\,\,i=k_1+k_2+1,\cdots,n.
\end{align*}
% We denote by $(\boldsymbol{b},\boldsymbol{a},\boldsymbol{L})$ all the reinsurance parameters, that is, \\
Let $(\boldsymbol{b},\boldsymbol{a},\boldsymbol{L})=(b_1,\cdots,b_{k_1},a_{k_1+1},\cdots,a_n,L_{k_1+k_2+1},\cdots,L_{n})$ denoted all the  reinsurance parameters.
% Suppose there there is a $\alpha_1$ such that
%$\boldsymbol{R}(\alpha_1)=x$
It is easy to see that 
$g(V,x,\boldsymbol{R})=g(V,x,\boldsymbol{b},\boldsymbol{a},\boldsymbol{L} )$ is a left-continuous function with negative jumps with respect to $b_i$'s and $a_i$'s, and right-continuous function with positive jumps  with respect to $L_i$'s. For example, let $k=n=3$, $p_{kk}=1$, $R_1(\alpha)=R(b_1,\alpha)= b_1\alpha$, $R_2(\alpha)=R(a_2,\alpha)= a_2\wedge\alpha$ and
$R_3(\alpha)=R(a_3,L_3,\alpha)= a_3\wedge\alpha+(\alpha-a_3-L_3)^+$ :
\begin{align*}
g(V,x,\boldsymbol{R})=& \frac{1}{\sum_{i=1}^{3} \beta_i}\bigg(\beta_1(1-p_{12})(1-p_{13})\int_{0}^{\infty}\big((1+\eta_1)\alpha +V(x-\alpha)\big)dp(R_1(U_1)\le \alpha)\\
&+\beta_2(1-p_{21})(1-p_{23})\int_{0}^{\infty}\big((1+\eta_1)\alpha +V(x-\alpha)\big)dp(R_2(U_2)\le \alpha)\\
&+{\beta_3} (1-p_{31})(1-p_{32})\int_{0}^{\infty}\big((1+\eta_1)\alpha +V(x-\alpha)\big)dp(R_3(U_3)\le \alpha)\\
&+
\big(\beta_1p_{12}(1-p_{13})+\beta_2p_{21}(1-p_{23})\big)\int_{0}^{\infty}\big((1+\eta_1)\alpha +V(x-\alpha)\big)dp(R_1(U_1)+R_2(U_2)\le \alpha)\\
&+
\big(\beta_1p_{13}(1-p_{12})+\beta_2p_{31}(1-p_{32})\big)\int_{0}^{\infty}\big((1+\eta_1)\alpha +V(x-\alpha)\big)dp(R_1(U_1)+R_3(U_3)\le \alpha)\\
&+
\big(\beta_2p_{23}(1-p_{21})+\beta_3p_{32}(1-p_{31})\big) \int_{0}^{\infty}\big((1+\eta_1)\alpha +V(x-\alpha)\big)dp(R_2(U_2)+R_3(U_3)\le \alpha)\\
&+\big(\beta_1p_{11}p_{12}p_{13}+\beta_2p_{21}p_{22}p_{23}+\beta_3 p_{31}p_{32}p_{33}\big)\\
&\times\int_{0}^{\infty}\big((1+\eta_1)\alpha +V(x-\alpha)\big)dp(R_1(U_1)+R_2(U_2)+R_3(U_3)\le \alpha)\bigg).
\end{align*}
We only investigate 
$$g_1(V,x,b_1,a_2,a_3,L_3)=\int_{0}^{\infty}\big((1+\eta_1)\alpha +V(x-\alpha)\big)dp(R_1(U_1)+R_2(U_2)+R_3(U_3)\le \alpha).$$
Suppose that there are $t_1$, $t_2$ and $t_3$ such that
 $R(b_1,t_1)+R(a_2,t_2)+R(a_3,L_3,t_3)=x$. It is easy to see that
 \begin{align*}
 R(b_1^-,t_1)+R(a_2,t_2)+R(a_3,L_3,t_3)&\le x \le R(b_1^+,t_1)+R(a_2,t_2)+R(a_3,L_3,t_3),\\
 R(b_1,t_1)+R(a_2^-,t_2)+R(a_3,L_3,t_3)&\le x \le R(b_1,t_1)+R(a_2^+,t_2)+R(a_3,L_3,t_3),\\
 R(b_1,t_1)+R(a_2,t_2)+R(a_3^-,L_3,t_3)&\le x \le R(b_1,t_1)+R(a_2,t_2)+R(a_3^+,L_3,t_3),\\
 R(b_1,t_1)+R(a_2,t_2)+R(a_3,L_3^+,t_3)&\le x \le R(b_1,t_1)+R(a_2,t_2)+R(a_3,L_3^-,t_3).
 \end{align*}
 So we obtain
   \begin{align*}
 g_1(V,x,b_1,a_2,a_3,L_3)=g_1(V,x,b_1^-,a_2,a_3,L_3)&\ge g_1(V,x,b_1^+,a_2,a_3,L_3),\\
  g_1(V,x,b_1,a_2,a_3,L_3)=g_1(V,x,b_1,a_2^-,a_3,L_3)&\ge g_1(V,x,b_1,a_2^+,a_3,L_3),\\
  g_1(V,x,b_1,a_2,a_3,L_3)=g_1(V,x,b_1,a_2,a_3^-,L_3)&\ge g_1(V,x,b_1,a_2,a_3^+,L_3),\\
 g_1(V,x,b_1,a_2,a_3,L_3)=g_1(V,x,b_1,a_2,a_3,L_3^+)&\le g_1(V,x,b_1^+,a_2,a_3,L_3^-),
  \end{align*}
   we have that $g_1(V,x,b_1,a_2,a_3,L_3)$ is a left-continuous function with negative jumps with respect to $b_1$, $a_2$ and $a_3$, and right-continuous function with positive jumps  with respect to $L_3$. So there exists  at least one vector $(\boldsymbol{b}^*,\boldsymbol{a}^*,\boldsymbol{L}^*)$ where the maximum of $g(V,x,\boldsymbol{b},\boldsymbol{a},\boldsymbol{L} )$ is attaind. $\square$

%%%%%%%%%%%%%%%%%%%%%%%%%%%%%%%%%%%%%%

{\bf  Proof of Proposition \ref{prop5.7}}\,\, we should  prove the following,
\begin{itemize}
	\item[(a)] 	${\cal{A}}^*$ is a closed set,
	\item[(b)]the lower limit of any connected component of ${\cal{B}}^{*}$ belongs to ${\cal{A}}^{*}$,
	\item[(c)] ${\cal{B}}^{*}$  is a left-open set,
	\item[(d)] there is an $y$ such that $(y , \infty) \subset {\cal{B}}^{*}$,
	\item[(e)] ${\cal{C}}^*$ is a right-open set, and
	\item[(f)]both ${\cal{A}}^*$ and ${\cal{B}}^{*}$ are nonempty.
\end{itemize}
 By Lemma \ref{lemma3.2.'}, we get that $${\cal{A}}^*=\{x\in \mathbb{R}_+ \,\,\text{suct that}\,\,\sup_{\boldsymbol{{R}} \in \boldsymbol{\mathcal{F}}}\,\,\Lambda_{\boldsymbol{R}}({V})(x)=0  \}$$ is closed.
 
 Now, we must show that  if $(x_1,x_2] \subset {\cal{B}}^{*}$ and $x_1\notin {\cal{B}}^{*}$ , then $x_1\in {\cal{A}}^{*}$, that is, the lower limit of any connected component of ${\cal{B}}^{*}$ belongs to ${\cal{A}}^{*}$.
 At first, we will take $x_1=0$. Let us define $Z_1$ as the severity of the first claim and $\tau_1$ as the time of this claim;
 then, consider the admissible strategy $\pi=(\boldsymbol{R}, D)\in \Pi_x$ as such that for $t\le \tau_1$, $D_t=p_{\boldsymbol{R}}t$. In this regard,
 \begin{eqnarray}\label{A.99}
 V_{\pi}(0)=E\left(\int_{0}^{\tau_1}e^{-cs}dD_s+e^{-c\tau_1}V_{\pi}\big(p_{\tau_1}-D_{\tau_1}-\boldsymbol{{R}}_{\tau_1}(Z_1)\big)\right)
 \end{eqnarray}
 Now, if $x$ belong to ${\cal{B}}^{*}$ and $x$ is immediately paid as a dividend, then $V'(x)=1$ and according to Lemma\ref{remark5.5} we have the following 
 $$V(0)=\lim_{x\rightarrow 0^+}\sup_{\pi\in \Pi_{0,x}} V_{\pi}(x)=\sup_{\pi\in \Pi_{0,0^+}} V_{\pi}(0).$$
 Moreover, $\Pi_{0,0^+}=\{\pi\in \Pi_{0} \,\,\text{such that }\,\, X_t^{\pi}<0^+\,\,\text{for all}\,\, t \ge 0\};$
 therefore the only possible state that remains is $D_t=p_{\boldsymbol{R}_0}t$ and here we obtain,
$$V(0)=\sup_{\boldsymbol{{R}}_0 }\frac{p_{\boldsymbol{R}_0}}{c+\sum_{i=1}^{m}\beta_i-\sum_{i=1}^{m}\beta_i p(\boldsymbol{R}_0(z)=0)}$$
Now, let consider 
$$\boldsymbol{R}_0^*=argmax \frac{p_{\boldsymbol{R}_0}}{c+\sum_{i=1}^{m}\beta_i-\sum_{i=1}^{m}\beta_i p(\boldsymbol{R}_0(z)=0)},$$ then $\Lambda(V)(0)=\Lambda_{\boldsymbol{R}_0^*}(V)(0)=0$ can be concluded; so, $0\in {\cal{A}}^{*}$.

Now, we will take $x_1\ge 0$. Here, on the one hand, it is crystal clear if $V'(x_1)=1$, then due to $x_1 \notin {\cal{B}}^{*}$, $x_1\in {\cal{A}}^{*}$ will be gleaned. On the other hand, let's consider $V'(x_1)\neq 1$. Here, we suppose 
$$\liminf_{x\rightarrow x_1^{-}}\frac{V(x)-V(x_1)}{x-x_1}=\delta >1,$$ then based upon the Definition \ref{Viscositydefi}, the following will be satisfied
$$ 
\max\left\{1-V'(x),\underset{\boldsymbol{R}}{sup}\,\,\left(p_{\boldsymbol{R}}\lambda-(c+\sum_{i=1}^m \beta_i)V(x)+(\sum_{i=1}^m \beta_i)\int_{0}^{x}V(x-\alpha)dG^{\boldsymbol{R}}(\alpha)\right)\right\}\ge 0.
$$
for all $\lambda \in (1,\delta]$.
So
$$\underset{\boldsymbol{R}}{sup}\,\,\left(p_{\boldsymbol{R}}\lambda-(c+\sum_{i=1}^m \beta_i)V(x_1)+(\sum_{i=1}^m \beta_i)\int_{0}^{x_1}V(x_1-\alpha)dG^{\boldsymbol{R}}(\alpha)\right)\ge 0,
$$
and therefor, 
$$\Lambda(x_1)=\underset{\boldsymbol{R}}{sup}\,\,\left(p_{\boldsymbol{R}}-(c+\sum_{i=1}^m \beta_i)V(x_1)+(\sum_{i=1}^m \beta_i)\int_{0}^{x_1}V(x_1-\alpha)dG^{\boldsymbol{R}}(\alpha)\right)\ge 0,
$$
for all $\lambda \in (1,\delta]$. Then, by Lemma \ref{lemma3.2.'}, we have that  $\Lambda(x_1)=0$.
Similar to the proof of Theorem 8.2 (c) by \cite{azcue2005optimal}, we can gain the similar result for the case 
$$\liminf_{x\rightarrow x_1^{-}}\frac{V(x)-V(x_1)}{x-x_1}=1.$$
Proof for cases (c), (d), (e) and (f) is similar to proof of Theorem 8.2 of \cite{azcue2005optimal}.$\square$

{\bf  Proof of Theorem \ref{theorem5.2}}\,\, 
Based on Lemma \ref{prop5.6} for $x \in  {\cal{A}}^*\bigcup{\cal{C}}^*$, there is a $\boldsymbol{R}_x^* \in \boldsymbol{\mathcal{F}}$ such that
$sup_{\boldsymbol{{R}} \in \boldsymbol{\mathcal{F}}}\,\,{\tilde{{\cal{L}}}}_{\boldsymbol{R}}({V})(x)={\tilde{{\cal{L}}}}_{\boldsymbol{R}_x^*}({V})(x).$ Now we define 
\begin{align*}
\tilde{\boldsymbol{R}}_x(\alpha)=
\begin{cases}	\boldsymbol{R}_x^*(\alpha)	& if \,\, x\in {\cal{A}}^*\bigcup{\cal{C}}^*\\
\boldsymbol{R}_0(\alpha)	& if \,\, x\in {\cal{B}}^*
\end{cases}
\end{align*}
Where $\boldsymbol{R}_0$ is any retained loss function. Furthermore, according to proposition \ref{prop5.7}, ${\cal{P}}^*$ is a partition. Now, 
		$\pi^*=\pi({\cal{P}}^*, \boldsymbol{R}^*)$ is demonstrated to be optimal, that is, $V(x)=V_{\pi^*}(x)$.
By using  Proposition \ref{prop5.5},  we should  prove the following,
\begin{itemize}
	\item[(a)] $V$ is left continuous at the upper limits of the connected component of ${\cal{C}}^*$,
	\item[(b)] $V$ is right continuous at the lower limits of the connected component of ${\cal{B}}^*$,
	\item[(c)] $V$ has derivative equal to one on ${\cal{B}}^*$,
	\item[(d)]$V$ is an almost-everywhere solution of ${\tilde{{\cal{L}}}}_{\boldsymbol{R}_x^*}({V})(x)$ in the connected components of ${\cal{C}}^*$,
	\item[(e)] $V$ is a solution of 
$$ \Lambda_{\boldsymbol{R}^*}(V)(x)=p_{\boldsymbol{R}^*}-(c+\sum_{i=1}^m \beta_i)V(x)+(\sum_{i=1}^m \beta_i)\int_{0}^{x}V(x-\alpha)dG^{\boldsymbol{R}^*}(\alpha),$$
in ${\cal{A}}^*$.
\end{itemize}
By Lemma \ref{lemma1ofcani}, $V$ is locally Lipschitz, so the cases (a) and (b) are true; by Definition \ref{defi5.3} and Lemma \ref{remark5.5}, $V'=1$ on ${\cal{B}}^*$ and $\Lambda_{\boldsymbol{R}^*}=0$ on ${\cal{A}}^*$, so (c) and (e) are true. Finally, $V$ is an almost-everywhere solution of ${\cal{L}}=0$ in the connected components of ${\cal{C}}^*$  because $V$ is a  viscosity solution of (\ref{hjbeqV}), and by Lemma \ref{lemma5.6}(b) $V'(x)>1$ at any $x$ of ${\cal{C}}^*$ where $V$ is differentiable,  hence (d) is true.  $\square$

  Before going in the proof of Theorem \ref{Theorem4.1}, we need the following two lemmas.
 \begin{lemma}\label{numeric-increasing}
	Let some small step size $ h $ such that $p_{{\boldsymbol{R}}}\ge 0, \,\, i=1,2, \cdots,$  and let $D_h= \{ih, i=0,\,1,\cdots\}$. Then 
	\  $f'_h(s) \ge 0$ for all $s\in D_h$.
\end{lemma}
{\bf Proof}
For $i=0$ the  assertion holds. 
Assume that $i$ is a positive integer with $f '_h(k h)\ge  0,  k = 1, ..., i- 1$. Then 
$f_h(0) \le f_h(h) \le \cdots \le f_h((i-1)h)$ and thus 
$$G_{\boldsymbol{R}}(ih)=\sum_{\{j\le i\}}f_h((i-j)h)P\{(j-1)h<\boldsymbol{R}(Y)\le jh\}\le f_h((i-1) h).$$
So for  $s=ih$
$$\frac{(\sum_{i=1}^n\beta_i) f_h(s-h)-(\sum_{i=1}^n\beta_i) G_{\boldsymbol{R}}(s)}{p_{{\boldsymbol{R}}}} \ge 0,$$
and thus obviously 
$$f'_h(ih)=f'_h(s)=\inf_{\boldsymbol{R}\in \mathcal{F}}\frac{(\sum_{i=1}^n\beta_i) f_h(s-h)-(\sum_{i=1}^n\beta_i) G_{\boldsymbol{R}}(s)}{p_{{\boldsymbol{R}}}} \ge 0,$$
which completes the induction.
$\square$

\begin{lemma}\label{numeric-visco}
	Let, ${\cal{A}}^h=\{a_1^h,a_2^h,\cdots, a_n^h\}$ and 
	${\cal{P}}=({\cal{A}},{\cal{B}},{\cal{C}})=\lim_{h\rightarrow 0}{\cal{P}}^h$.
	In the setting of the above lemma, and define
	\begin{align}
	V^{*}(s)=\lim_{h\rightarrow 0}\sup_{ih \rightarrow s} V_h(ih),
	\end{align}
	and 
	\begin{align}\label{numeric_supersolution}
	V_{*}(s)=\lim_{h\rightarrow 0}\inf_{ih \rightarrow s} V_h(ih),
	\end{align}
	if $V_h(x)\le x+c$, then the functions $V^{*}(s)$ and $V_{*}(s)$ are respectively, sub and super viscosity solution of (\ref{hjbeqV}).
\end{lemma}
{\bf proof}  Firstly, we show that the  function $V^*$ is locally Lipschitz and a viscosity subsolution of \ref{hjbeqV}. Fix $M>0$ and let $x_1$ and $x_2$ that are belong to the ${\cal{A}}\cup {\cal{C}}$ and $0\le x_1<x_2<M$. Take sequences $k_{n}^{(i)}h_n $, $i=1,2$, such that
$$k_{n}^{(i)}h_n \rightarrow x_i, \,\,i=1,2,\qquad V(k_{n}^{(i)}h_n) \rightarrow V(x_i) $$
It is easy to
see that
\begin{align*}
V^*(x_2)-V^*(x_1)=& \lim_{n\longrightarrow \infty} V_{h_n}(k_n^{(2)}h_n)-V^*(x_1)\\
& \le \lim_{n\longrightarrow \infty} V_{h_n}(k_n^{(2)}h_n)- \lim_{n\longrightarrow \infty}  V_{h_n}(k_n^{(1)}h_n)\\
& \le \lim_{n\longrightarrow \infty} \big(V_{h_n}(k_n^{(2)}h_n)-V_{h_n}(k_n^{(1)}h_n)\big)\\
&\le \limsup_{n\rightarrow \infty}\frac{V_{h_n}(k_n^{(2)} h_n)-V_{h_n}(k_n^{(1)}h_n)}{k_n^{(2)} h_n-k_n^{(1)}h_n}|k_n^{(1)} h_n-k_n^{(2)}h_n|\\
&=\limsup_{n\rightarrow \infty}\frac{h_n\sum_{k_n^{(1)} h_n < ih_n \le k_n^{(2)}h_n }V_{h_n}'(i h_n)}{k_n^{(2)} h_n-k_n^{(1)}h_n}|k_n^{(2)} h_n-k_n^{(1)}h_n|\\
&=\limsup_{n\rightarrow \infty}\frac{\sum_{k_n^{(1)} h_n < ih_n \le k_n^{(2)}h_n }V_{h_n}'(i h_n)}{k_n^{(2)} -k_n^{(1)}}|k_n^{(2)} h_n-k_n^{(1)}h_n|.
\end{align*}
But we have 
\begin{align*}
f'_{ih_n}(ih_n)
\le \frac{(\delta+(\sum_{i=1}^n\beta_i)) (f_{h}(ih_n-h_n)-(\sum_{i=1}^n\beta_i) G(ih_n))}{p}
\le \frac{(\delta+(\sum_{i=1}^n\beta_i)) f_{h}(ih_n-h_n)}{p}
\end{align*}
Now, according to  the above inequality and $V_{h_n}(ih_n)\le ih_n+c$, the following relation is obtained:
$$ V^*(x_2)-V^*(x_1) \le K(x_2-x_1),$$
where
$$K\le\frac{\delta+\sum_{i=1}^n\beta_i}{\sum_{i=1}^n\beta_i}(M+c).$$
For the case that $x_1$ or $x_2$, or both do not belong to ${\cal{A}}\cup {\cal{C}}$, proof is obvious.

To show that $V^{*}$ is a viscosity subsolution,  suppose that $w$ is a test function such that $ V^*(x)-w(x)$ has a maximum at $s > 0$.  Note that, for $h$ sufficiently small we can find $h<kh\in D_h$ such that
$ V_{h}(kh+h)- V_{h}(kh) \le w(kh+h) - w(kh).$ So, we have
\begin{align} &\underset{\boldsymbol{\mathcal{F}}}{sup}\,\,\{p_{{\boldsymbol{R}}}\frac{w(kh+h) - w(kh)}{h}-(\sum_{i=1}^n\beta_i+\delta) V_{h}(kh+h)+(\sum_{i=1}^n\beta_i) G_{\boldsymbol{R}}(kh+h)\}\nonumber\\
&\ge  \underset{\boldsymbol{\mathcal{F}}}{sup}\,\,\{p_{{\boldsymbol{R}}}\frac{V_{h}(kh+h)- V_{h}(kh)}{h}-(\sum_{i=1}^n\beta_i+\delta) V_{{h}}(kh+h)+(\sum_{i=1}^n\beta_i) G_{\boldsymbol{R}}(kh+h)\}\nonumber\\
&\ge  
p_{{\boldsymbol{R}}_h}\frac{V_{h}(kh+h)- V_{h}(kh)}{h}-(\sum_{i=1}^n\beta_i+\delta) V_{{h}}(kh+h)+(\sum_{i=1}^n\beta_i) G_{\boldsymbol{R}_h}(kh+h)\nonumber\\\
&= 0.
\label{AA1}
\end{align}
Take  sequences $k_n$ and $h_n$ such that $h_n\longrightarrow 0$, $V_{h_{n}}(k_{n}h_{n})\longrightarrow V^*(s)$  and $k_{n}h_{n}\longrightarrow s$.
Then by Fatou’s lemma,
\begin{align}
\limsup_{n\Rightarrow \infty} G_{\boldsymbol{R}}(k_{n}h_{n}+h_{n}) \le E_{\boldsymbol{R}}(V^*(s-Y)), \qquad V_{{h_{n}}}(k_{n}h_{n}+h_{n})\rightarrow V^*(k_{n}h_{n}).
\label{AA2}
\end{align}
So, from \ref{AA1} and \ref{AA2}, we have
\begin{align*}
0\le&  \underset{\boldsymbol{\mathcal{F}}}{sup}\,\,\{p_{{\boldsymbol{R}}}w'(k_{n}h_{n})-(\sum_{i=1}^n\beta_i+\delta) V_{{h_{n}}}(k_{n}h_{n}+h_{n})(k_{n}h_{n})+(\sum_{i=1}^n\beta_i) G_{\boldsymbol{R}}(k_{n}h_{n}+h_{n})\}\\
\le & \underset{\boldsymbol{\mathcal{F}}}{sup}\,\,\{p_{{\boldsymbol{R}}}w'(k_{n}h_{n})-(\sum_{i=1}^n\beta_i+\delta) V^*(k_{n}h_{n})+(\sum_{i=1}^n\beta_i) E_{\boldsymbol{R}}(V^*(x-Y))\}\\
\le & \underset{\boldsymbol{\mathcal{F}}}{sup}\,\,\{p_{{\boldsymbol{R}}}w'(s)-(\sum_{i=1}^n\beta_i+\delta)w(s)+(\sum_{i=1}^n\beta_i) E_{\boldsymbol{R}}(w(s-Y))\}.
\end{align*}
Thus, $V^*$ is a viscosity subsolution. Similarly, $V_*$  is locally Lipschitz and a viscosity supersolution of \ref{hjbeqV}.
$\square$

{\bf Proof of Theorem \ref{Theorem4.1}.}
By Lemma \ref{numeric-visco}, $V^*$ and $V_*$ are locally Lipschitz and
based on Lemma \ref{numeric-increasing} and the numerical algorithm described above, $V_*$ and $V^*$ are greater than one, and therefore satisfy the conditions (i) and (ii)  of Definition \ref{growthA1}.
Moreover, according to the condition  $V_h(x)\le{x}+{c}$, $V^*$ and $V_*$ are bounded from above by ${x}+{c}$. Hence, $V^*$ and $V_*$  belong to  $L$.
In the proof of Proposition \ref{prop4.4} we show that $V(.)$  is smaller or equal than any supersolution of (\ref{hjbeqV}) that belongs to $L$. So, $V \le V_*$.	
On the other hand, given that 	${\cal{P}}=\lim_{h\rightarrow 0}{\cal{P}}^h$  is a band partition and $V^*$ is a Lipschitz function, it is to see that  $V^*$ is satisfied in the conditions of Proposition \ref{prop5.5}, and therefore $V^*=V_{{\cal{P}}}$, where $\pi=\pi({\cal{P}},\boldsymbol{R})$ is
reinsurance band strategy associated with ${\cal{P}}$  and $\boldsymbol{R}$.
Therefore
$V^*\le V$.
Since  $V^*\ge V_*$ by definition, we have convergence and therefore $V(s)= V^*(s)$ $(=V_*(s))$. 
$\square$
\bibliographystyle{agsm} 
\bibliography{references1.bib}
\end{document}